\documentclass{conm-p-l}

\usepackage{amsmath,amssymb,latexsym}
\usepackage{booktabs,tabularx}
\usepackage{graphicx,psfrag}
\usepackage{url}
\usepackage{txfonts}

\newtheorem{deff}{Definition}
\newtheorem{lem}[deff]{Lemma}
\newtheorem{prop}[deff]{Proposition}
\newtheorem{thm}[deff]{Theorem}

\newtheorem{conj}[deff]{Conjecture}

\makeatletter

\begin{document}

\title{Graph Coloring Manifolds}

\author{P\'eter Csorba}
\address{Department of Mathematics,
The University of Western Ontario,
London, Ontario N6A 5B7, Canada}
\email{pcsorba@uwo.ca}

\author{Frank H.~Lutz}
\address{Technische Universit\"at Berlin,
Institut f\"ur Mathematik, 
Str.\ des 17.\ Juni 136,
10623 Berlin, Germany}
\email{lutz@math.tu-berlin.de}

\thanks{The first author was supported by the joint Berlin/Z\"{u}rich graduate
         program ``Combinatorics, Geometry, and Computation'',
         by grants from NSERC and the Canada Research Chairs program.}



\begin{abstract}
We introduce a new and rich class of 
\emph{graph coloring manifolds}
via the Hom complex construction of Lov\'asz.
The class comprises examples of Stiefel manifolds,
series of spheres and products of spheres, cubical surfaces,
as well as examples of Seifert manifolds. Asymptotically,
graph coloring manifolds provide examples of highly connected, 
highly symmetric manifolds. 
\end{abstract}

\maketitle

\section{Introduction}

In the topological approach to graph coloring, initiated by Lov\'asz'  proof \cite{Lovasz1978} 
of the Kneser Conjecture~\cite{MKneser1955}, lower bounds on the \emph{chromatic number} $\chi(H)$ 
of a graph $H$ are obtained by exploiting topolo\-gical invariants of
a simplicial or cell complex $K(H)$ that is associated with~$H$.

There are several standard constructions that associate a topological 
space $K(H)$ with a graph~$H$, e.g., the (simplicial) neighborhood complex ${\mathcal N}(H)$
of Lov\'asz \cite{Lovasz1978}, the (simplicial) box complex $B(H)$ 
of Matou\v{s}ek and Ziegler \cite{MatousekZiegler2004},
and, with respect to a reference graph $G$, 
the (cellular) Hom complex ${\rm Hom}(G,H)$ of Lov\'asz 
(cf.\ \cite{BabsonKozlov2006}, \cite{Kozlov2005cpre}).

From an algorithmic point of view, the topological approach seems,
up to now, not suitable to produce ``good'' lower bounds
on $\chi(H)$ for general input graphs $H$:
For example, the historically first topological lower bound by Lov\'asz 
requires the computation of the connectivity of the
neighborhood complex ${\mathcal N}(H)$. 
\begin{thm} {\rm (Lov\'asz~\cite{Lovasz1978})}
Let $H$ be a graph. If ${\mathcal N}(H)$ is $k$-connected, then $\chi(H)\geq k+3$.
\end{thm}
However, neighborhood complexes of graphs can be of
``arbitrary'' homotopy type~\cite{Csorba2004pre},
and for general complexes it is \emph{not decidable} 
whether they are $1$-connected or not! 
Moreover, there are cases for which the
connectivity could be determined, but for which
the corresponding lower bounds are far from tight
\cite{Walker1983}.

It is therefore most surprising that for highly
structured, highly symmetric graphs such as Kneser graphs
and generalization \cite{AlonFranklLovasz1986} the topological approach provides sharp
lower bounds while other approaches fail badly; \cite{MatousekZiegler2004}
discusses this issue and gives further references.

In order to get away from connectivity, lower bounds have been
formulated in terms of topological invariants that \emph{are}
computable \cite[Remark~2.7]{BabsonKozlov2004pre}, 
or the topological tools have been replaced by 
purely combinatorial ones; see Matou\v{s}ek~\cite{Matousek2004}. 
Still, the size of the associated complexes causes problems, 
since for the standard constructions the number of cells 
of the complexes $K(H)$ grows exponentially.

In a recent series of papers, Babson and Kozlov \cite{BabsonKozlov2006}, \cite{BabsonKozlov2004pre}
and Kozlov~\cite{Kozlov2005cpre} intensively studied properties of Hom complexes ${\rm Hom}(G,H)$
and proved new topological lower \linebreak
bounds (see as well \v{C}uki\'c and Kozlov \cite{CukicKozlov2004pre},
\cite{CukicKozlov2005}, Schultz \cite{Schultz2005pre}, and \v{Z}ivaljevi\'c \cite{Zivaljevic2005pre}). 
For example, it turned out that ${\rm Hom}(K_2,K_n)$ is a PL sphere of dimension $n-2$, for $n\geq 2$,
and by spectral sequence calculations that the Hom complexes ${\rm Hom}(C_5,K_n)$
have the (co)homology of Stiefel manifolds. This was the starting
point for the first author to formulate Conjecture~\ref{conj:csorba}
(see Section~\ref{sec:one}) that the Hom complexes ${\rm Hom}(C_5,K_n)$ \emph{are} (PL) \emph{homeomorphic} 
to Stiefel manifolds.

In this paper, we will show that the Hom complexes ${\rm Hom}(C_5,K_n)$
indeed are PL manifolds. More generally, we will characterize in
Theorem~\ref{thm:main} (Section~\ref{sec:main})
those graphs $G$ for which the Hom complexes ${\rm Hom}(G,K_n)$ are PL manifolds for all
$n\geq\chi(G)$. Such manifolds we call \emph{graph coloring manifolds}.

In Section~\ref{sec:basic}, we give a short account on Hom complexes. 
Section~\ref{sec:main} introduces graph coloring manifolds.
Various examples and series of examples of graph coloring manifolds
are discussed in Sections~\ref{sec:zero}--\ref{sec:two}.

\section{Basic Definitions, Notations, and Examples}
\label{sec:basic}

Let $G=(V(G),E(G))$ and $H=(V(H),E(H))$ be two graphs with node sets
$V(G)$ and $V(H)$ and edge sets $E(G) \subseteq \binom{V(G)}{2}$
and $E(H)\subseteq \binom{V(H)}{2}$, respectively.
We assume that the graphs are \emph{simple graphs}, 
i.e., graphs without loops and parallel edges. 

A \emph{graph homomorphism} is a map $\phi :V(G)\rightarrow V(H)$,
such that if $\{i,j\}\in E(G)$, 
then $\{\phi(i),\phi(j)\}\in E(H)$, 
that is, the image of every edge of the graph $G$ is an edge 
of the graph~$H$. Let the set of all graph homomorphisms 
from $G$ to $H$ be denoted by ${\mathcal H}om(G,H)$. 
For two disjoint sets of vertices $A,B \subseteq V(G)$ we define
$G[A,B]$ as the 
subgraph of $G$ with $V(G[A,B])= A \cup B$ and
$E(G[A,B])=\{ \{a,b\} \in E(G) \colon a\in A,b\in B\}$.

Let $\Delta^{V(H)}$ be the (abstract) simplex whose set of vertices is $V(H)$.
Furthermore, let $C(G,H)$ denote the direct product $\prod_{x\in V(G)}\Delta^{V(H)}$,
i.e., the copies of $\Delta^{V(H)}$ are indexed by vertices of $G$.
A cell of $C(G,H)$ is a direct product of simplices $\prod_{x\in V(G)}\sigma_x$.

\begin{deff}
For any pair of graphs $G$ and $H$ let the
\emph{Hom complex} ${\rm Hom}(G,H)$
be a subcomplex of $C(G,H)$ defined by the following
condition: $c=\prod_{x\in V(G)}\sigma_x\in{\rm Hom}(G,H)$ if and only if
for any $u,v\in V(G)$ if $\{u,v\}\in E(G)$, then $H[\sigma_u,\sigma_v]$
is complete bipartite.
\end{deff}

\noindent
The topology of ${\rm Hom}(G,H)$ is inherited from the product
topology of $C(G,H)$.
Thus, ${\rm Hom}(G,H)$ is a polyhedral complex whose
(non-empty) cells are products of simplices and are indexed by functions
(multi-homomorphisms)
$\eta :V(G)\rightarrow 2^{V(H)}\backslash\{\emptyset\}$,
such that if $\{\imath,\jmath\}\in E(G)$, then for every
$\tilde{\imath}\in\eta (\imath)$ and $\tilde{\jmath}\in\eta (\jmath)$
it follows that $\{\tilde{\imath},\tilde{\jmath}\}\in E(H)$.

Let $V(G)=\{ 1,\dots,m\}$. We encode the functions $\eta$ 
by vectors $(\eta (1),\dots,\eta (m))$ of non-empty sets 
$\eta (i)\subseteq V(H)$ with the above properties.
A cell $(A_1,\dots,A_m)$ of ${\rm Hom}(G,H)$ is a \emph{face}
of a cell $(B_1,\dots,B_m)$ of ${\rm Hom}(G,H)$
if $A_i\subseteq B_i$ for all $1\leq i\leq m$.
In particular, ${\rm Hom}(G,H)$ has ${\mathcal H}om(G,H)$ as its set of vertices.
Moreover, every cell $(A_1,\dots,A_m)$ of ${\rm Hom}(G,H)$ 
is a product of $m$ simplices of dimension $|A_i|-1$ for $1\leq i\leq m$.
For brevity, we write sets $A=\{ a_1,\dots,a_k\} \subseteq V(H)$
in compressed form as strings, i.e., $A=a_1\dots a_k$.

A cell of ${\rm Hom}(G,H)$ is a \emph{maximal face} or \emph{facet} 
if it is not contained in any higher-dimensional cell of ${\rm Hom}(G,H)$.

\medskip
\smallskip

\begin{figure}
\begin{center}
\footnotesize
\psfrag{(1,2)}{(1,2)}
\psfrag{(1,3)}{(1,3)}
\psfrag{(2,3)}{(2,3)}
\psfrag{(2,1)}{(2,1)}
\psfrag{(3,1)}{(3,1)}
\psfrag{(3,2)}{(3,2)}
\psfrag{(12,3)}{(12,3)}
\psfrag{(1,23)}{(1,23)}
\psfrag{(13,2)}{(13,2)}
\psfrag{(3,12)}{(3,12)}
\psfrag{(23,1)}{(23,1)}
\psfrag{(2,13)}{(2,13)}
\includegraphics[width=.475\linewidth]{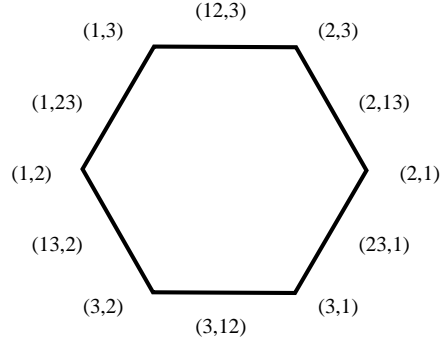}
\end{center}
\caption{The Hom complex ${\rm Hom}(K_2,K_3)$.}
\label{fig:Hom_K2_K3}
\end{figure}

\noindent
\emph{Example 1:}\, The cells of the Hom complex ${\rm Hom}(K_2,K_3)$
are given by the vectors $(1,2)$, $(1,3)$, $(2,3)$, $(2,1)$, $(3,1)$, $(3,2)$, 
$(12,3)$, $(13,2)$, $(23,1)$, $(3,12)$, $(2,13)$, and $(1,23)$.
Therefore, ${\rm Hom}(K_2,K_3)$ is a cycle with six edges;
see Figure~\ref{fig:Hom_K2_K3}.

\bigskip

\noindent
\emph{Example 2:}\, 
The Hom complex ${\rm Hom}(K_2,K_n)$ is a PL sphere of dimension $n{-}2$
for $n\geq 2$. In fact, ${\rm Hom}(K_2,K_n)$ is the boundary complex of
a polytope \cite[Sect.~4.2]{BabsonKozlov2006}: it can be described
as the boundary of the Minkowski sum of an $({n-1})$-dimensional simplex
$\sigma_{n-1}$ and its negative $-\sigma_{n-1}$, as stated in
\cite[p.~107, Ex.~3 (c)]{Matousek2003}.

\section{Vertex-Links and Flag Simplicial Spheres}
\label{sec:main}

Babson and Kozlov asked in \cite{BabsonKozlov2006}
for what graphs the Hom complex construction provides 
a connection to polytopes. In this section, 
we will characterize those graphs $G$ for which 
${\rm Hom}(G,K_n)$ is a piecewise linear (PL) manifold
for all $n\geq\chi(G)$.

A (finite) simplicial complex is a \emph{PL $d$-manifold}
if and only if every vertex-link is a PL $(d-1)$-sphere,
i.e., every vertex-link is PL homeomorphic to the boundary
of the standard $d$-simplex $\sigma_d$.


There are several ways to define the \emph{link} of a vertex $v$ for
polyhedral complexes. For Hom complexes ${\rm Hom}(G,H)$ we will use the
following. Let the face poset of ${\rm Hom}(G,H)$ be denoted by $\mathcal{F}({\rm Hom}(G,H))$
and let the link of $v$ in ${\rm Hom}(G,H)$ be the cell complex whose face poset is given by
$\mathcal{F}_{>v}({\rm Hom}(G,H))$. This link then is a simplicial
complex since ${\rm Hom}(G,H)$ is a prodsimplicial complex (cf.\ \cite[2.4.3]{Kozlov2005cpre}).

For a graph $G$ we say that $X\subseteq V(G)$ is an \emph{independent set}
if there is no edge between any two vertices of $X$.
The \emph{independent set complex} ${\rm Ind}(G)$ of a graph $G$ is
the simplicial complex with vertex set $V(G)$ and $X\subseteq V(G)$
forming a simplex if and only if X is an independent set in $G$, i.e.,
${\rm Ind}(G)=\{X\subseteq V(G)\,|\,X \mbox{ is independent in $G$}\}$.

Every cell of ${\rm Hom}(G,K_n)$ corresponds to a multi-coloring
$f\colon V(G) \to 2^{\{1,\dots,n\}}\backslash\{\emptyset\}$, 
where the map $f$ assigns $f(v)$ distinct colors to every vertex $v\in V(G)$,
such that the set of vertices colored by any color $i\in\{1,\dots,n\}$ 
forms an independent set in $G$. We denote these sets by
$\Delta_i(f)= \{v\in V(G)\,|\,i\in f(v)\}$ and consider
them as simplices of ${\rm Ind}(G)$.

\begin{lem}
Let $\phi$ be a vertex of ${\rm Hom}(G,K_n)$, 
i.e., a proper coloring of $G$ which we regard 
as a multi-coloring  $\phi\colon V(G) \to 2^{\{1,\dots,n\}}\backslash\{\emptyset\}$
with $|\phi(v)|=1$ for all $v\in V(G)$.
Then
\[
{\rm link}_{\,{\rm Hom}(G,K_n)}(\phi)
\]
is isomorphic to the join product
\[
{\rm link}_{\,{\rm Ind}(G)}(\Delta_1(\phi)) * \cdots * {\rm link}_{\,{\rm Ind}(G)}(\Delta_n(\phi)).
\]

\end{lem}

\noindent
\textbf{Proof.}
A simplex of the first complex ${\rm link}_{\,{\rm Hom}(G,K_n)}(\phi)$
corresponds to a multi-coloring $f\colon V(G) \to 2^{\{1,\dots,n\}}\backslash\{\emptyset\}$ 
which extends $\phi$.
We can consider such an extension color-wise. 
For color $i\in\{1,\dots,n\}$ we have that
$\Delta_i(\phi)\subseteq\Delta_i(f)\in{\rm Ind}(G)$
and therefore $\Delta_i(f)\backslash\Delta_i(\phi)\in{\rm link}_{\,{\rm Ind}(G)}(\Delta_i(\phi))$.
Thus we can identify $f$ with
$(\Delta_1(f)\backslash\Delta_1(\phi),\dots,\Delta_n(f)\backslash\Delta_n(\phi))$
and therefore can regard $f$ as an element of
${\rm link}_{\,{\rm Ind}(G)}(\Delta_1(\phi)) * \cdots * {\rm link}_{\,{\rm Ind}(G)}(\Delta_n(\phi))$.
Conversely, every simplex of ${\rm link}_{\,{\rm Ind}(G)}(\Delta_1(\phi)) * \cdots * {\rm link}_{\,{\rm Ind}(G)}(\Delta_n(\phi))$
gives rise to a unique extension $f$ of $\phi$.
\mbox{}\hfill$\Box$

\begin{lem}\label{lem:referee}
\mbox{}\\[-5mm]
\begin{enumerate}
\item If\, ${\rm Ind}(G)$ is a PL sphere, then ${\rm Hom}(G,K_n)$ is a PL manifold for any $n\geq \chi(G)$.
\item If\, ${\rm Hom}(G,K_n)$ is a PL manifold and $n>\chi(G)$, then ${\rm Ind}(G)$ is a PL sphere.
\end{enumerate}
\end{lem}

\noindent
\textbf{Proof.}
1. Let ${\rm Ind}(G)$ be a PL sphere. Since the link of any simplex of
a PL sphere is a PL sphere (of lower dimension)
and since the join product of PL spheres is again a PL sphere,
it follows by the previous lemma that the link 
of any vertex of ${\rm Hom}(G,K_n)$ is a PL sphere.
Thus, ${\rm Hom}(G,K_n)$ is a PL manifold.

2. Let ${\rm Hom}(G,K_n)$ be a PL manifold.
Since $n>\chi(G)$, there is a vertex $\phi$ of ${\rm Hom}(G,K_n)$
that does not use the color $n$. Hence,
${\rm link}_{{\rm Ind}(G)}(\Delta_n(\phi))={\rm Ind}(G)$.
Since ${\rm Hom}(G,K_n)$ is a PL manifold,
${\rm link}_{\,{\rm Hom}(G,K_n)}(\phi) \cong {\rm link}_{\,{\rm Ind}(G)}(\Delta_1(\phi)) * \cdots * {\rm link}_{\,{\rm Ind}(G)}(\Delta_n(\phi))$
is a PL sphere. Now, the join product of simplicial complexes is a PL sphere
if and only if every factor is a PL sphere (see \cite[2.24(5)]{RourkeSanderson1982}).
It follows that the last factor, ${\rm link}_{{\rm Ind}(G)}(\Delta_n(\phi))={\rm Ind}(G)$,
is a PL sphere.\hfill$\Box$

\medskip
We can formulate this result in terms of $G$ using the following definition.

\begin{deff}
Let $K$ be a (finite) simplicial complex. If $K$ has no ``empty simplices'', i.e., 
if every set of vertices of $K$ which form a clique in the
$1$-skeleton ${\rm Skel}_1(K)$ actually spans a simplex, then $K$ is a \emph{flag simplicial complex}
(cf.\ \cite{CharneyDavis1995}).
A \emph{flag simplicial sphere} is a flag simplicial complex 
which triangulates a sphere.
\end{deff}

The \emph{clique complex} ${\rm{Cliq}}(G)=\{X\subseteq V(G)\,|\,X \mbox{ is a clique in $G$}\}$
of any graph $G$ is a flag simplicial complex in a natural way
with $G={\rm Skel}_1({\rm{Cliq}}(G))$.

\begin{thm}\label{thm:main}
Let $G$ be a graph. Then the Hom complex\, ${\rm Hom}(G,K_n)$ 
is a PL manifold for all $n\geq \chi(G)$ 
if and only if $G$ is the complement of the $1$-skeleton 
of a flag simplicial PL sphere.
\end{thm}

\noindent
\textbf{Proof.}
Let ${\rm Hom}(G,K_n)$ be a PL manifold for all $n\geq \chi(G)$.
Then, in particular, \linebreak
 ${\rm Hom}(G,K_{\chi(G)+1})$ is a PL manifold,
and thus, by Lemma~\ref{lem:referee}, ${\rm Ind}(G)={\rm Cliq}(\overline{G})$ is a PL sphere.
Hence, $G$ is the complement of the $1$-skeleton of the flag simplicial
PL sphere~${\rm Cliq}(\overline{G})$.

Conversely, if $G$ is the complement of the $1$-skeleton of a flag
simplicial PL sphere $K$, i.e., $G=\overline{{\rm Skel}_1(K)}$,
then ${\rm Ind}(G)={\rm Cliq}({\rm Skel}_1(K))=K$ is a flag simplicial
PL sphere and therefore ${\rm Hom}(G,K_n)$ a PL manifold by Lemma~\ref{lem:referee}.\hfill$\Box$

\medskip

\noindent
\emph{Remark 1:}\,
If $n<\chi(G)$, then ${{\rm Hom}}(G,K_n)=\emptyset$.
If $n=\chi(G)$, then every vertex $\phi$ of ${{\rm Hom}}(G,K_{\chi(G)})$
uses all colors $1,\dots,\chi(G)$. If ${{\rm Hom}}(G,K_{\chi(G)})$
is a PL manifold, then ${\rm{Ind}}(G)$ need not be
a PL sphere. It is only required, that 
the links of vertices (or of higher-dimensional faces if 
every color is used more than once in every vertex of
${{\rm Hom}}(G,K_{\chi(G)})$) of ${\rm{Ind}}(G)$
are flag simplicial PL spheres.
In particular, if $G$ is the complement of the 1-skeleton
of a flag combinatorial manifold, then ${{\rm Hom}}(G,K_{\chi(G)})$
is a PL manifold. As another example, if $G$ is a connected bipartite graph, then
${{\rm Hom}}(G,K_2)=S^0$.

\medskip
\smallskip

\noindent
\emph{Remark 2:}\, It is possible for ${{\rm Hom}}(G,K_n)$ to be a (non-PL) manifold,
even without ${\rm{Ind}}(G)$ being a sphere. (See \cite{BjoernerLutz2000} for a discussion
of non-PL spheres and non-PL manifolds.) For example,
if ${\rm{Ind}}(G)$ is a flag combinatorial homology sphere
(i.e., a combinatorial manifold with the homology of a sphere, but 
not homeomorphic to the standard sphere)
and $n>\chi(G)$, then for every vertex $\phi$ of ${\rm Hom}(G,K_n)$
the join product 
${\rm link}_{\,{\rm Hom}(G,K_n)}(\phi) \cong {\rm link}_{\,{\rm Ind}(G)}(\Delta_1(\phi)) * \cdots * {\rm link}_{\,{\rm Ind}(G)}(\Delta_n(\phi))$
is a simplicial sphere by the double
suspension theorem of Edwards \cite{Edwards1975} and Cannon
\cite{Cannon1979}. Also, if $G$ is the complement of the 1-skeleton
of a flag simplicial non-PL sphere, then ${{\rm Hom}}(G,K_n)$ is a non-PL manifold
for $n>\chi(G)$.


\begin{deff}\label{deff:graphcoloringmanifold}
A Hom complex ${\rm Hom}(G,K_n)$ is a \emph{graph coloring manifold}
if $G$ is the complement of the $1$-skeleton of a flag simplicial PL sphere.
\end{deff}

\noindent
\emph{Remark 3:}\, By Definition~\ref{deff:graphcoloringmanifold}
and Theorem~\ref{thm:main} graph coloring manifold are PL manifolds.

\medskip
\smallskip

\noindent
\emph{Remark 4:}\, Graph coloring manifolds are highly symmetric:
relabeling the colors of $K_n$ defines an action of the 
symmetric group $S_n$ on ${\rm Hom}(G,K_n)$.

\medskip
\smallskip

Babson and Kozlov \cite[2.4]{BabsonKozlov2006} stated as a basic
property of Hom complexes that
\begin{equation}\label{eq:direct_product}
{{\rm Hom}}(G_1\,\Dot{\cup}\,\,G_2,H)={{\rm Hom}}(G_1,H)\times {{\rm Hom}}(G_2,H),
\end{equation}
from which it follows that if\,
$G=\Dot{\displaystyle\bigcup_{i=1,\dots,k}}K_2$\, is the complement
of the $1$-skeleton of the boundary of the $k$-dimensional crosspolytope $\partial C_k^{\Delta}$,
then 
\begin{equation}
{{\rm Hom}}\,(\Dot{\displaystyle\bigcup_{i=1,\dots,k}}K_2,K_n)=\varprod\limits_{i=1,\dots,k}S^{n-2}.
\end{equation}

\begin{deff}
A flag simplicial PL sphere is \emph{prime} if the complement 
of its $1$-skeleton is connected. A Hom complex ${\rm Hom}(G,K_n)$
is a \emph{graph coloring manifold of sphere dimension~$d$} 
if $G$ is the complement of the $1$-skeleton of a prime flag
simplicial PL sphere of dimension~$d$.
\end{deff}

Since every coloring of a graph $G$ can be regarded as a covering of $G$ 
by independent sets, the following lower bound holds
for the chromatic number $\chi(G)$ of $G$:
\begin{equation}
\chi(G)\geq\Bigl\lceil\frac{|V|}{\alpha(G)}\Bigl\rceil=\Bigl\lceil\frac{|V|}{\omega(\overline{G})}\Bigl\rceil,
\end{equation}
where $\alpha(G)$ is the \emph{independence number} or \emph{stable set number} of $G$
(i.e., the maximum size of an independent set in~$G$)
and $\omega(G)$ is the \emph{clique number} of $G$
(i.e., the maximum size of a clique in $G$).
If $G$ is the complement of the $1$-skeleton of a prime flag
simplicial $PL$ $d$-sphere
on $m$ vertices, then\, $\alpha(G)=\omega(\overline{G})=d+1$. Thus
\begin{equation}
\label{eqn:flag}
\chi(G)\geq \Bigl\lceil\frac{m}{d+1}\Bigl\rceil
\end{equation}
and
\begin{equation}
{\rm dim}({\rm Hom}(G,K_{\chi(G)+k}))=(\chi(G)+k)(d+1)-m
\end{equation}
for all\, $k\geq 0$.

The lower bound~(\ref{eqn:flag}) can be arbitrarily bad: 
If $G$ is the complement of the $1$-skeleton of the
suspension $S^0*C_{2r+1}$ of an odd cycle $C_{2r+1}$, $r\geq 2$, 
then\, $\chi(G)=2r+1>\lceil\frac{2r+3}{3}\rceil$.

From the following theorem it follows that graph coloring
manifolds provide examples of highly connected manifolds.
\begin{thm}  {\rm (\v{C}uki\'c and Kozlov~\cite{CukicKozlov2005})}
Let $G$ be a graph of maximal valency $s$, then the Hom complex ${{\rm Hom}}(G,K_n)$
is at least $(n-s-2)$-connected.
\end{thm}
Let $G$ be the complement of the 1-skeleton
of a flag simplicial PL sphere.
If $G$ has maximal valency $s$, then ${{\rm Hom}}(G,K_n)$
is simply connected and thus orientable for $n\geq s+3$.
We expect that ${{\rm Hom}}(G,K_n)$ is orientable also for\, $\chi(G)\leq n< s+3$.

\begin{conj}
Graph coloring manifolds are orientable.
\end{conj}

\section{Graph Coloring Manifolds of Sphere Dimension Zero}
\label{sec:zero}

Trivially, $S^0$, consisting of two isolated vertices, 
is the only zero-dimensional flag simplicial sphere.
The complement of its (empty) 1-skeleton is the complete graph $K_2$.
Hence, the graph coloring manifolds of sphere dimension zero
are the Hom complexes ${\rm Hom}(K_2,K_n)\cong S^{n-2}$, for $n\geq 2$.

\section{Graph Coloring Manifolds of Sphere Dimension One}
\label{sec:one}

The one-dimensional flag simplicial spheres
are the cycles $C_m$ of length $m\geq 4$. 
For $m=4$ we have that (the $1$-skeleton) ${\rm SK}_1(\overline{C}_4)=\overline{C}_4=K_2\,\Dot{\cup}\,K_2$
with
$${\rm Hom}(K_2\,\Dot{\cup}\,K_2,K_n)={\rm Hom}(K_2,K_n)\times{\rm Hom}(K_2,K_n)\cong S^{n-2}\!\times S^{n-2}.$$
If $m\geq 5$, then ${\rm SK}_1(\overline{C}_m)=\overline{C}_m$ is connected.
In the following, we treat odd and even cycles separately.

\subsection{Hom Complexes of Complements of Odd Cycles}

Babson and Kozlov \cite{BabsonKozlov2004pre} used topological
information on the Hom complexes ${\rm Hom}(C_5,K_n)$
(with $\overline{C}_5\cong C_5$ for $m=5$) 
and, more generally, on the Hom complexes ${\rm Hom}(C_{2r+1},K_n)$, for $r\geq 2$ and $n\geq r+1$, 
to prove the Lov\'asz Conjecture:
\begin{thm} {\rm (Babson and Kozlov \cite{BabsonKozlov2004pre})}
If for a graph $H$ the complex\, ${\rm Hom}(C_{2r+1},H)$ is $k$-connected, 
for some $r\geq 1$ and $k\geq -1$, then $\chi(H)\geq k+4$.
\end{thm}
Babson and Kozlov computed various cohomology groups of the Hom complexes ${\rm Hom}(C_m,K_n)$.
For $m=5$, the respective cohomology groups are those of Stiefel manifolds.
\begin{conj} {\rm (Csorba)}\label{conj:csorba}
The Hom complex ${{\rm Hom}}(C_5,K_{n+2})$ 
is PL homeomorphic to the Stiefel mani\-fold $V_{n+1,2}$.
\end{conj}
It is elementary to verify that ${{\rm Hom}}(C_5,K_3)$
consists of two cycles with $15$ vertices and $15$ edges each.

\medskip
\smallskip

\noindent
\emph{Example 5:}\, ${{\rm Hom}}(C_5,K_3)\cong V_{2,2}\cong S^0\times S^1$.

\medskip
\smallskip

For $n=2$, the complex ${{\rm Hom}}(C_5,K_4)$ has $240$ vertices
and $300$ maximal cells that are either cubes or prisms
over triangles.

Since every cell of a Hom complex is a product of simplices,
triangulations of graph coloring manifolds (\emph{without additional vertices}) 
can easily be obtained
by the product triangulation construction as described in \cite{Lutz2003bpre}.
For small examples, the homology of the resulting triangulations
can then be computed with one of the programs \cite{DumasHeckenbachSaundersWelker2003b}
or~\cite{polymake}. 

The product triangulation of ${\rm Hom}(C_5,K_4)$  has 
$f$-vector $f=(240,1680,2880,1440)$. As homology we obtained
$H_*({\rm Hom}(C_5,K_4))=({\mathbb Z},{\mathbb Z}_2,0,{\mathbb Z})$,
which coincides with the spectral sequence computations
of Babson and Kozlov in \cite{BabsonKozlov2004pre}.
We also computed the homology of ${\rm Hom}(C_7,K_4)$
and obtained that $H_*({\rm Hom}(C_7,K_4))=({\mathbb Z},{\mathbb Z}_2,0,{\mathbb Z},0,0)$
which was conjectural in \cite{BabsonKozlov2004pre}.

We next used the bistellar flip heuristic BISTELLAR \cite{Lutz_BISTELLAR}
to determine that the complex ${\rm Hom}(C_5,K_4)$ is homeomorphic to ${\mathbb R}{\bf P}^3$.
(See \cite{BjoernerLutz2000} for a discussion of the heuristic;
for large complexes the \texttt{bistellar} client (due to N. Witte) 
of the TOPAZ module of the \texttt{polymake} system~\cite{polymake} provides 
a fast implementation of BISTELLAR.)

\begin{thm}\label{thm:rp3}
\mbox{}\,${{\rm Hom}}(C_5,K_4)\cong V_{3,2}\cong {\mathbb R}{\bf P}^3$. 
\end{thm}

\noindent
\textbf{Proof.}
In addition to the above computer proof, we give an explicit Heegaard
decomposition of ${{\rm Hom}}(C_5,K_4)$ from which one can see that 
this Hom complex is homeomorphic to~${\mathbb R}{\bf P}^3$
(and thus homeomorphic to the Stiefel manifold $V_{3,2}$).

First we show that the collection of cells of the form $(ijk,*,*,*,*)$
forms a solid torus. By symmetry it is enough to consider the
collection of cells $(123,*,*,*,*)$. Since the numbers 1, 2, and 3 can not occur at positions 2 and 5,
it immediately follows that the cells of this collection are of the form $(123,4,*,*,4)$.
The middle $(*,*)$-part is the six-gon corresponding to ${\rm Hom}(K_2,K_3)$;
see Figure~\ref{fig:Hom_K2_K3}. So $(ijk,*,*,*,*)$ is the product of a triangle and a circle.

\begin{figure}
\begin{center}
\footnotesize
\psfrag{(124,3,24,1,3)}{\tiny (124,3,24,1,3)}
\psfrag{(124,3,4,12,3)}{\tiny (124,3,4,12,3)}
\psfrag{(124,3,14,2,3)}{\tiny (124,3,14,2,3)}
\psfrag{(124,3,1,24,3)}{\tiny (124,3,1,24,3)}
\psfrag{(124,3,12,4,3)}{\tiny (124,3,12,4,3)}
\psfrag{(124,3,2,14,3)}{\tiny (124,3,2,14,3)}
\psfrag{(123,4,2,13,4)}{\tiny (123,4,2,13,4)}
\psfrag{(123,4,12,3,4)}{\tiny (123,4,12,3,4)}
\psfrag{(123,4,1,23,4)}{\tiny (123,4,1,23,4)}
\psfrag{(123,4,13,2,4)}{\tiny (123,4,13,2,4)}
\psfrag{(123,4,3,12,4)}{\tiny (123,4,3,12,4)}
\psfrag{(123,4,23,1,4)}{\tiny (123,4,23,1,4)}
\psfrag{(12,34,2,1,34)}{\tiny (12,34,2,1,34)}
\psfrag{(12,34,1,2,34)}{\tiny (12,34,1,2,34)}
\psfrag{(12,3,24,1,34)}{\tiny (12,3,24,1,34)}
\psfrag{(12,3,14,2,34)}{\tiny (12,3,14,2,34)}
\psfrag{(12,34,1,24,3)}{\tiny (12,34,1,24,3)}
\psfrag{(12,34,2,14,3)}{\tiny (12,34,2,14,3)}
\psfrag{(12,34,2,13,4)}{\tiny (12,34,2,13,4)}
\psfrag{(12,34,1,23,4)}{\tiny (12,34,1,23,4)}
\psfrag{(12,4,13,2,34)}{\tiny (12,4,13,2,34)}
\psfrag{(12,4,23,1,34)}{\tiny (12,4,23,1,34)}

\psfrag{(12,3,4,12,34)}{\tiny (12,3,4,12,34)}
\psfrag{(12,3,4,123,4)}{\tiny (12,3,4,123,4)}
\psfrag{(12,3,24,13,4)}{\tiny (12,3,24,13,4)}
\psfrag{(12,3,14,23,4)}{\tiny (12,3,14,23,4)}
\psfrag{(12,3,124,3,4)}{\tiny (12,3,124,3,4)}
\psfrag{(12,34,12,3,4)}{\tiny (12,34,12,3,4)}
\psfrag{(12,34,12,4,3)}{\tiny (12,34,12,4,3)}
\psfrag{(12,4,123,4,3)}{\tiny (12,4,123,4,3)}
\psfrag{(12,4,13,24,3)}{\tiny (12,4,13,24,3)}
\psfrag{(12,4,23,14,3)}{\tiny (12,4,23,14,3)}
\psfrag{(12,4,3,124,3)}{\tiny (12,4,3,124,3)}
\psfrag{(12,4,3,12,34)}{\tiny (12,4,3,12,34)}
\includegraphics[width=.95\linewidth]{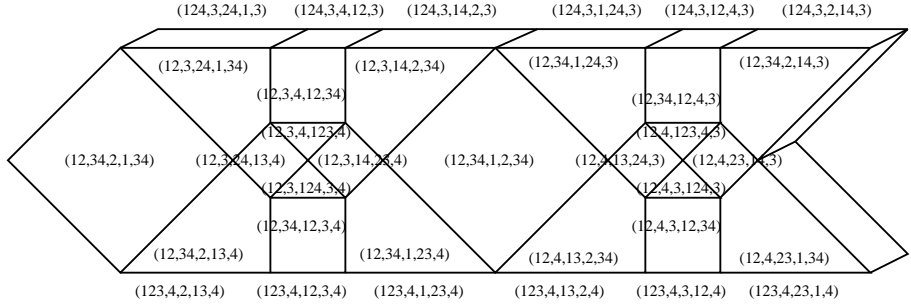}
\end{center}
\caption{The solid torus $(12,*,*,*,*)$ in ${\rm Hom}(C_5,K_4)$.}
\label{fig:collection_12}
\end{figure}

\begin{figure}
\begin{center}
\footnotesize
\psfrag{(12,34,2,13,4)}{\tiny (1,34,2,13,4)}
\psfrag{(12,3,24,1,34)}{\tiny (1,3,24,1,34)}
\psfrag{(12,3,4,12,34)}{\tiny (1,3,4,12,34)}
\psfrag{(12,34,12,3,4)}{\tiny (1,34,12,3,4)}
\psfrag{(12,3,4,123,4)}{\tiny (1,3,4,123,4)}
\psfrag{(12,3,124,3,4)}{\tiny (1,3,124,3,4)}
\psfrag{(12,3,24,13,4)}{\tiny (1,3,24,13,4)}
\psfrag{(12,3,14,2,34)}{\tiny (1,3,14,2,34)}
\psfrag{(12,3,14,23,4)}{\tiny (1,3,14,23,4)}
\psfrag{(12,34,1,23,4)}{\tiny (1,34,1,23,4)}
\psfrag{(12,34,2,14,3)}{\tiny (1,34,2,14,3)}
\psfrag{(12,4,23,1,34)}{\tiny (1,4,23,1,34)}
\psfrag{(12,34,12,4,3)}{\tiny (1,34,12,4,3)}
\psfrag{(12,4,3,12,34)}{\tiny (1,4,3,12,34)}
\psfrag{(12,4,23,14,3)}{\tiny (1,4,23,14,3)}
\psfrag{(12,4,123,4,3)}{\tiny (1,4,123,4,3)}
\psfrag{(12,4,3,124,3)}{\tiny (1,4,3,124,3)}
\psfrag{(12,4,13,24,3)}{\tiny (1,4,13,24,3)}
\psfrag{(12,34,1,24,3)}{\tiny (1,34,1,24,3)}
\psfrag{(12,4,13,2,34)}{\tiny (1,4,13,2,34)}
\psfrag{(12,34,1,2,34)}{\tiny (1,34,1,2,34)}
\psfrag{(12,34,2,1,34)}{\tiny (1,34,2,1,34)}

\psfrag{(13,24,13,2,4)}{\tiny (1,24,13,2,4)}
\psfrag{(13,24,1,23,4)}{\tiny (1,24,1,23,4)}
\psfrag{(13,4,12,3,24)}{\tiny (1,4,12,3,24)}
\psfrag{(13,4,2,13,42)}{\tiny (1,4,2,13,24)}
\psfrag{(13,4,23,1,24)}{\tiny (1,4,23,1,24)}
\psfrag{(13,2,134,2,4)}{\tiny (1,2,134,2,4)}
\psfrag{(13,2,4,123,4)}{\tiny (1,2,4,123,4)}
\psfrag{(13,2,34,1,24)}{\tiny (1,2,34,1,24)}
\psfrag{(13,2,34,12,4)}{\tiny (1,2,34,12,4)}
\psfrag{(13,2,14,23,4)}{\tiny (1,2,14,23,4)}
\psfrag{(13,2,4,13,24)}{\tiny (1,2,4,13,24)}
\psfrag{(13,2,14,3,24)}{\tiny (1,2,14,3,24)}
\psfrag{(13,24,1,3,24)}{\tiny (1,24,1,3,24)}
\psfrag{(13,24,1,34,2)}{\tiny (1,24,1,34,2)}
\psfrag{(13,24,13,4,2)}{\tiny (1,24,13,4,2)}
\psfrag{(13,24,3,14,2)}{\tiny (1,24,3,14,2)}
\psfrag{(13,4,2,134,2)}{\tiny (1,4,2,134,2)}
\psfrag{(13,4,12,34,2)}{\tiny (1,4,12,34,2)}
\psfrag{(13,4,23,14,2)}{\tiny (1,4,23,14,2)}
\psfrag{(13,4,123,4,2)}{\tiny (1,4,123,4,2)}
\psfrag{(13,24,3,1,24)}{\tiny (1,24,3,1,24)}
\psfrag{(13,24,3,12,4)}{\tiny (1,24,3,12,4)}

\psfrag{(14,3,24,1,23)}{\tiny (1,3,24,1,23)}
\psfrag{(14,3,24,13,2)}{\tiny (1,3,24,13,2)}
\psfrag{(14,3,2,134,2)}{\tiny (1,3,2,134,2)}
\psfrag{(14,3,124,3,2)}{\tiny (1,3,124,3,2)}
\psfrag{(14,3,12,34,2)}{\tiny (1,3,12,34,2)}
\psfrag{(14,23,1,34,2)}{\tiny (1,23,1,34,2)}
\psfrag{(14,23,1,4,23)}{\tiny (1,23,1,4,23)}
\psfrag{(14,2,13,4,23)}{\tiny (1,2,13,4,23)}
\psfrag{(14,23,1,24,3)}{\tiny (1,23,1,24,3)}
\psfrag{(14,2,34,12,3)}{\tiny (1,2,34,12,3)}
\psfrag{(14,2,13,24,3)}{\tiny (1,2,13,24,3)}
\psfrag{(14,2,3,124,3)}{\tiny (1,2,3,124,3)}
\psfrag{(14,2,134,2,3)}{\tiny (1,2,134,2,3)}
\psfrag{(14,23,14,2,3)}{\tiny (1,23,14,2,3)}
\psfrag{(14,2,3,14,23)}{\tiny (1,2,3,14,23)}
\psfrag{(14,23,4,13,2)}{\tiny (1,23,4,13,2)}
\psfrag{(14,3,12,4,32)}{\tiny (1,3,12,4,23)}
\psfrag{(14,3,2,14,23)}{\tiny (1,3,2,14,23)}
\psfrag{(14,23,14,3,2)}{\tiny (1,23,14,3,2)}
\psfrag{(14,23,4,1,23)}{\tiny (1,23,4,1,23)}
\psfrag{(14,32,4,12,3)}{\tiny (1,23,4,12,3)}
\psfrag{(14,2,34,1,23)}{\tiny (1,2,34,1,23)}

\includegraphics[width=.95\linewidth]{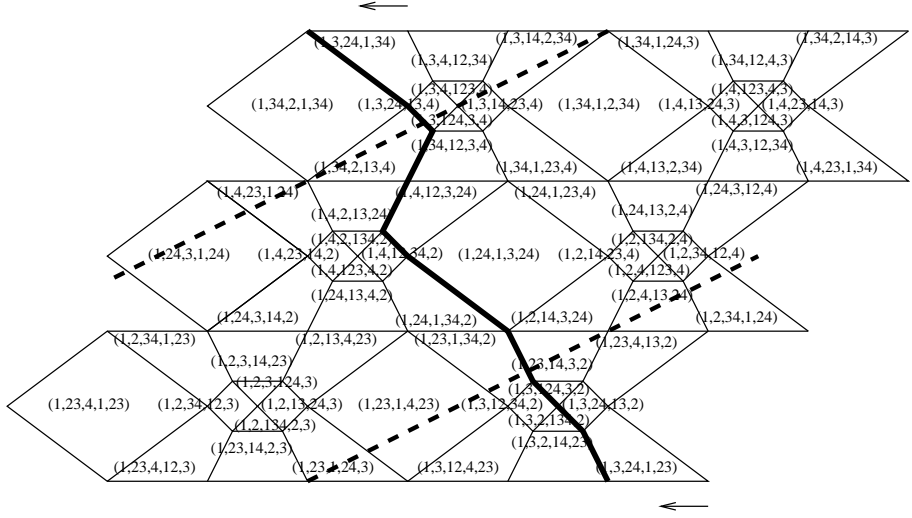}
\end{center}
\vspace{-0.3cm}
\caption{The boundary of the solid torus $(1,*,*,*,*)$ in ${\rm Hom}(C_5,K_4)$.}
\label{fig:cell1xxxx}
\end{figure}

In Figure~\ref{fig:collection_12} we display the collections of cells of the form
$(12,*,*,*,*)$ (with the cell on the left glued to the cells on the right).
Clearly, this collection of cells forms a torus, and therefore, by
symmetry, also every collection $(ij,*,*,*,*)$.
Finally, the cells of the form $(i,*,*,*,*)$ form a solid torus as well.
The boundary torus of $(1,*,*,*,*)$ can be seen in Figure~\ref{fig:cell1xxxx}.
Again, the left side is glued to the right side of Figure~\ref{fig:cell1xxxx}.
The gluing of the top and bottom is indicated by the arrows.

\begin{figure}
\begin{center}
\footnotesize
\psfrag{(1,3,2,1,3)}{\footnotesize (1,3,2,1,3)}
\psfrag{(1,2,1,3,2)}{\footnotesize (1,2,1,3,2)}
\psfrag{(1,3,2,3,2)}{\footnotesize (1,3,2,3,2)}
\psfrag{(1,3,2,1,2)}{\footnotesize (1,3,2,1,2)}
\psfrag{(1,3,2,1,4)}{\footnotesize (1,3,2,1,4)}
\psfrag{(1,3,2,3,4)}{\footnotesize (1,3,2,3,4)}
\psfrag{(1,3,2,1,234)}{\footnotesize (1,3,2,1,234)}
\psfrag{(1,4,2,3,4)}{\footnotesize (1,4,2,3,4)}
\psfrag{(1,4,2,3,2)}{\footnotesize (1,4,2,3,2)}
\psfrag{(1,4,1,3,2)}{\footnotesize (1,4,1,3,2)}
\psfrag{(1,234,1,3,2)}{\footnotesize (1,234,1,3,2)}
\psfrag{(1,34,12,3,2)}{\footnotesize (1,34,12,3,2)}
\psfrag{(1,34,2,3,24)}{\footnotesize (1,34,2,3,24)}
\psfrag{(1,3,2,13,24)}{\footnotesize (1,3,2,13,24)}
\psfrag{(1,3,1,3,2)}{\footnotesize (1,3,1,3,2)}
\includegraphics[width=.65\linewidth]{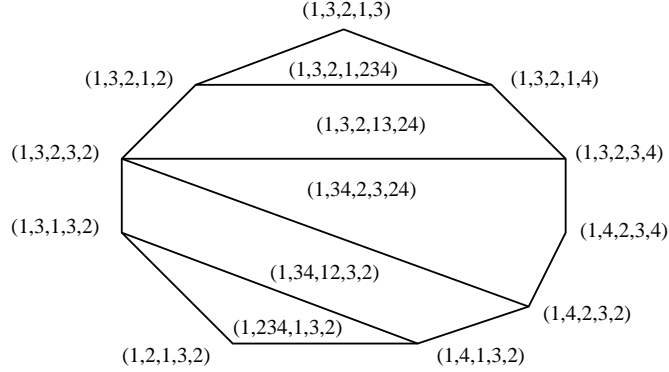}
\end{center}
\vspace{-0.3cm}
\caption{A meridian disk of the solid torus $(1,*,*,*,*)$ in ${\rm Hom}(C_5,K_4)$.}
\label{fig:meridian}
\end{figure}

\begin{figure}
\begin{center}
\footnotesize
\psfrag{123}{123}
\psfrag{124}{124}
\psfrag{134}{134}
\psfrag{234}{234}
\psfrag{12}{12}
\psfrag{13}{13}
\psfrag{14}{14}
\psfrag{23}{23}
\psfrag{24}{24}
\psfrag{34}{34}
\psfrag{1}{1}
\psfrag{2}{2}
\psfrag{3}{3}
\psfrag{4}{4}
\psfrag{23==(23,*,*,*,*)  etc.}{}
\includegraphics[width=.6\linewidth]{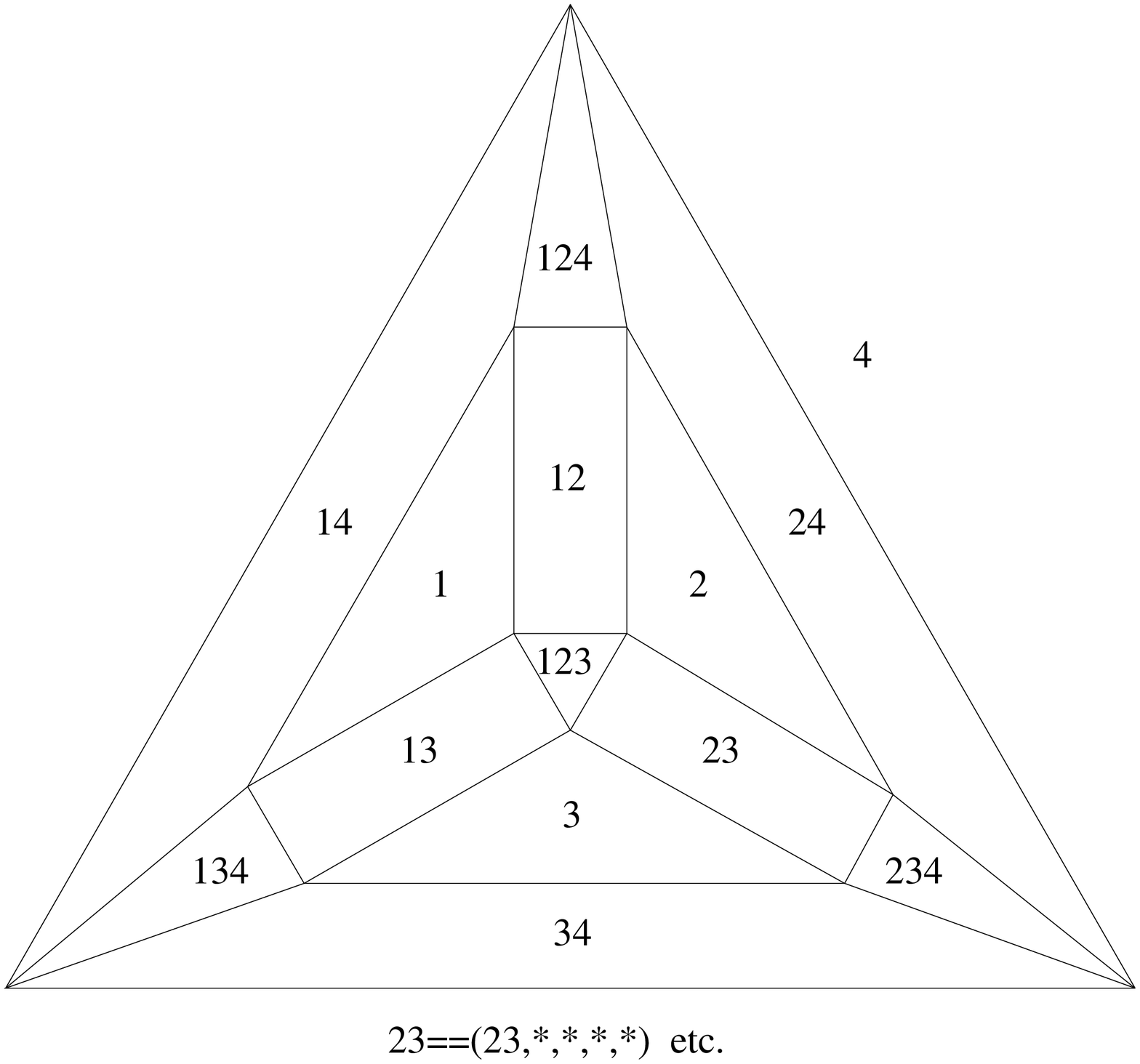}
\end{center}
\caption{Cell decomposition of the base sphere $S^2$ in ${\rm Hom}(C_5,K_4)$.}
\label{fig:base_S2}
\end{figure}

In order to understand how these solid tori are glued together we have
to identify meridian disks. For the collections $(ijk,*,*,*,*)$ and $(ij,*,*,*,*)$ 
this is clear. A meridian disk of the collection $(1,*,*,*,*)$ is
given in  Figure~\ref{fig:meridian}; its boundary corresponds to the
thick line in Figure~\ref{fig:cell1xxxx}.
\noindent
The complement of $(1,*,*,*,*)$ in ${\rm Hom}(C_5,K_4)$ is a solid torus
composed of the collections $(12,*,*,*,*)$, $(13,*,*,*,*)$, \dots, $(2,*,*,*,*)$, \dots, $(234,*,*,*,*)$,
which we abbreviate by $12$, $13$, \dots, $2$, \dots, $234$ in Figure~\ref{fig:base_S2}.
In fact, Figure~\ref{fig:base_S2} gives the base sphere of the
$S^1$-fibered space ${\rm Hom}(C_5,K_4)$ and makes clear how the
different tori are glued together.

A meridian curve of the complement of $(1,*,*,*,*)$
is drawn as a dashed curve in Figure~\ref{fig:cell1xxxx}.
Since this curve is a $(2,1)$-curve, it follows that
${\rm Hom}(C_5,K_4)$ is homeomorphic to the $3$-dimensional
real projective space. The latter space
is homeomorphic to the Stiefel manifold $V_{3,2}$.\mbox{}\hfill$\Box$

\medskip

The $5$-dimensional Hom complex ${{\rm Hom}}(C_5,K_5)$ consists of
$2070$ maximal cells and has $1020$ vertices. 
The corresponding product triangulation
has $f$-vector $(1020,25770,143900,$ 
$307950,283200,94400)$.
With the \texttt{bistellar} client it took less than a week to
reduce this triangulation to a triangulation with
$f=(12,66,220,390,336,112)$. 
The latter triangulation is $3$-neighborly, i.e., it has a complete
$2$-skeleton, and thus is simply connected. Its homology is
$({\mathbb Z},0,{\mathbb Z},{\mathbb Z},0,{\mathbb Z})$.
Moreover, its second Stiefel-Whitney class is trivial,
as we computed with \texttt{polymake}.
By the classification of simply connected $5$-manifolds
of Barden \cite{Barden1965}, the unique simply connected $5$-manifold
with homology $({\mathbb Z},0,{\mathbb Z},{\mathbb Z},0,{\mathbb Z})$
and trivial second Stiefel-Whitney class is $S^3\times S^2$.

\begin{thm}
\mbox{}\,${{\rm Hom}}(C_5,K_5)\cong V_{4,2}\cong S^3\times S^2$. 
\end{thm}

In the following, we discuss a particular representation of odd cycles
that gives some insight into \emph{all} Hom complexes ${\rm Hom}(\overline{C}_{2r+1},K_n)$
of complements of odd cycles $\overline{C}_{2r+1}$, $r\geq 2$.
(With a similar approach we will analyze the Hom complexes ${\rm Hom}(\overline{C}_{2r},K_n)$ 
of complements of even cycles $\overline{C}_{2r}$, $r\geq 2$, in the next section.)

\begin{figure}
\begin{center}
\footnotesize
\psfrag{a}{$a$}
\psfrag{b}{$b$}
\psfrag{c}{$c$}
\psfrag{A}{$A$}
\psfrag{B}{$B$}
\psfrag{C}{$C$}
\psfrag{D}{$D$}
\includegraphics[height=4.25cm]{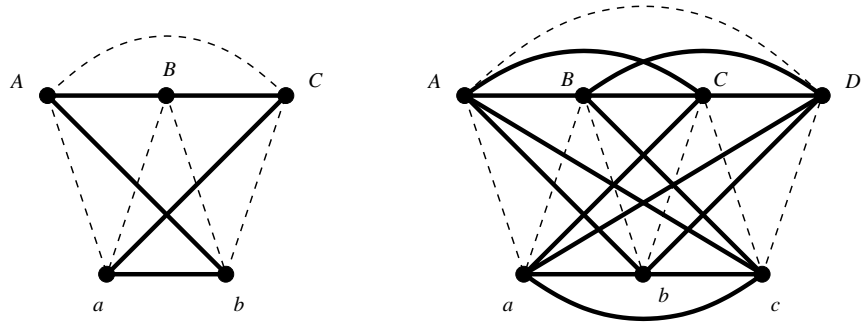}
\end{center}
\caption{The (dashed) cycles $C_5$ and $C_7$ and their complements.}
\label{fig:crown}
\end{figure}

We display the cycles $C_{2r+1}$, $r\geq 2$, in form of a crown that
is turned upside down; see Figure~\ref{fig:crown} for the \emph{crown representations} 
of the (dashed) cycles $C_5$ and $C_7$. Clearly, the bottom vertices of a crown representation 
form a clique, i.e., a complete graph $K_r$, 
in the complement $\overline{C}_{2r+1}$. 
 
Let us have a look at the crown representation of $C_5$.
Every cell $(a,b,A,B,C)$ of ${\rm Hom}(\overline{C}_5,K_n)$
contains every number $x\in\{1,\dots,n\}$ at exactly two
positions. Since the sets $a$ and $b$ are associated
with the bottom vertices that form a clique $K_2$ in $\overline{C}_5$,
the number $x$ can appear in at most one of the sets $a$ and $b$.
If it is contained in, say, $a$, then the second copy of $x$
can only be placed in the sets $A$ and $B$ that are connected
with $a$ by a dashed edge of $C_5$. The top vertices of $\overline{C}_5$
form a clique minus the (dashed) edge between the leftmost vertex 
and the rightmost vertex. Hence, if $x$ is contained in neither $a$
nor $b$, then it is contained in the leftmost top set $A$ and
in the rightmost top set $C$.

\begin{figure}
\begin{center}
\footnotesize
\psfrag{(12,3)}{$(12,3)$}
\psfrag{(13,2)}{$(13,2)$}
\psfrag{(23,1)}{$(23,1)$}
\psfrag{(3,12)}{$(3,12)$}
\psfrag{(2,13)}{$(2,13)$}
\psfrag{(1,23)}{$(1,23)$}
\psfrag{(1,2)}{$(1,2)$}
\psfrag{(1,3)}{$(1,3)$}
\psfrag{(2,3)}{$(2,3)$}
\psfrag{(2,1)}{$(2,1)$}
\psfrag{(3,1)}{$(3,1)$}
\psfrag{(3,2)}{$(3,2)$}
\includegraphics[height=4.5cm]{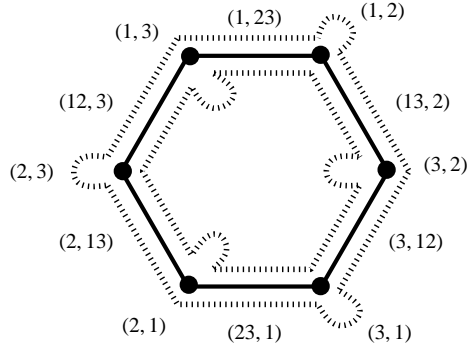}
\end{center}
\caption{${\rm Hom}(K_2,K_3)$ and ${\rm Hom}(\overline{C}_5,K_3)$.}
\label{fig:C_5_crown}
\end{figure}

If we restrict us further to $n=3$ colors,
then ${\rm Hom}(K_2,K_3)$ is a six-gon 
as displayed in solid in Figure~\ref{fig:C_5_crown}.
The cell $(a,b)=(1,23)$ of ${\rm Hom}(K_2,K_3)$
can be extended to a cell $(a,b;A,B,C)$ of ${\rm Hom}(\overline{C}_5,K_3)$
in precisely two ways, either to $(1,23;1,2,3)$ or to $(1,23;1,3,2)$.
We depict these edges of ${\rm Hom}(\overline{C}_5,K_3)$
as dashed edges in Figure~\ref{fig:C_5_crown}, parallel 
to the edge $(1,23)$ of ${\rm Hom}(K_2,K_3)$.
Let $(1,23;1,2,3)$ be the upper dashed edge. If we move the number $3$
from the second to the third position, then we obtain
the cell $(1,2;13,2,3)$ from which we move on to $(1,2;3,12,3)$,
and from there to $(1,2;3,1,23)$. These three cells of ${\rm Hom}(\overline{C}_5,K_3)$
correspond to the vertex $(1,2)$ of ${\rm Hom}(K_2,K_3)$
and are displayed together by a dashed half-cycle at the vertex $(1,2)$ 
in Figure~\ref{fig:C_5_crown}. If we move on further, then we get
to the dashed edge $(13,2;3,1,2)$, from there to the dashed
edge $(3,12;3,1,2)$, before we again start a half-cycle
$(3,1;23,1,2)$, $(3,1;2,13,2)$, $(3,1;2,3,12)$,
this time at the vertex $(3,1)$ of ${\rm Hom}(K_2,K_3)$.
We can then continue on the outer dashed cycle until 
we reach our starting edge $(1,23;1,2,3)$ of ${\rm Hom}(\overline{C}_5,K_3)$.
Similarly, we can move around the inner dashed cycle
when we start with $(1,23;1,3,2)$.

\begin{prop}
The Hom complex ${\rm Hom}(\overline{C}_{2r+1},K_{r+1})$
is the disjoint union of $r!$ cycles with  
$(2r^2+3r+1)$ vertices each.
\end{prop}

\noindent
\textbf{Proof.}
We first count the number of vertices of ${\rm Hom}(\overline{C}_{2r+1},K_{r+1})$,
i.e., the number of distinct colorings with $r+1$ colors of $\overline{C}_{2r+1}$.
To color the bottom $K_r$ in the crown representation of $\overline{C}_{2r+1}$
we choose $r$ of the $r+1$ colors and then have $r!$ choices
to place these $r$ colors. For one such coloring, say $(1,2,\dots,r)$,
there are $(2r+1)$ ways to extend it to a coloring of $\overline{C}_{2r+1}$:
If we use the color $r+1$ just once, then we have $r+1$ choices 
to place it in the top row of the crown; the remaining positions for
the colors in the top row are then completely determined by the position of the color
$(r+1)$ and by our choice of the colors in the bottom row.
If we use the color $r+1$ twice, then we have to put it at the positions
$1$ and $r+1$ of the top row. We further choose one of the colors $1,\dots,r$
not to be used in the top row; this again determines all the positions
for the colors in the top row. Thus we have $(r+1)$ choices if color $r+1$ appears
once in the top row and $r$ choices if color $r+1$ appears
twice in the top row. Altogether we have
$$\binom{r+1}{r}r!(r+1+r)=(2r^2+3r+1)r!$$
different colorings of $\overline{C}_{2r+1}$ with $r+1$ colors.

Since every number $1,\dots,r+1$ appears exactly twice
in a cell of ${\rm Hom}(\overline{C}_{2r+1},K_{r+1})$,
the dimension of ${\rm Hom}(\overline{C}_{2r+1},K_{r+1})$
is $2(r+1)-(2r+1)=1$. 
If we move for the edge $(1,2,\dots,r-1,r(r+1);1,2,\dots,r,r+1)$
of ${\rm Hom}(\overline{C}_{2r+1},K_{r+1})$
the number $r+1$ from the last position of the bottom
row to the first position of the top row and then continue
until we reach the edge $(r+1,1,2,\dots,(r-1)r;r+1,1,2,\dots,r-1,r)$,
this takes $r+1+r=2r+1$ steps. After $r+1$ rounds
we return to the starting edge $(1,2,\dots,r-1,r(r+1);1,2,\dots,r-1,r,r+1)$.
Thus, by symmetry, every cycle of ${\rm Hom}(\overline{C}_{2r+1},K_{r+1})$
has length $(2r+1)(r+1)=2r^2+3r+1$. Since ${\rm Hom}(\overline{C}_{2r+1},K_{r+1})$
has $(2r^2+3r+1)r!$ vertices, it follows that ${\rm Hom}(\overline{C}_{2r+1},K_{r+1})$
consists of $r!$ cycles with $(2r^2+3r+1)$ vertices each.\hfill$\Box$

\begin{figure}
\begin{center}
\footnotesize
\psfrag{(1,2,3)}{(1,2,3)}
\psfrag{(1,2,4)}{(1,2,4)}
\psfrag{(1,3,4)}{(1,3,4)}
\psfrag{(1,3,2)}{(1,3,2)}
\psfrag{(1,4,2)}{(1,4,2)}
\psfrag{(1,4,3)}{(1,4,3)}
\psfrag{(2,1,3)}{(2,1,3)}
\psfrag{(2,1,4)}{(2,1,4)}
\psfrag{(3,1,4)}{(3,1,4)}
\psfrag{(3,1,2)}{(3,1,2)}
\psfrag{(4,1,2)}{(4,1,2)}
\psfrag{(4,1,3)}{(4,1,3)}
\psfrag{(2,3,1)}{(2,3,1)}
\psfrag{(2,4,1)}{(2,4,1)}
\psfrag{(3,4,1)}{(3,4,1)}
\psfrag{(3,2,1)}{(3,2,1)}
\psfrag{(4,2,1)}{(4,2,1)}
\psfrag{(4,3,1)}{(4,3,1)}
\psfrag{(2,3,4)}{(2,3,4)}
\psfrag{(2,4,3)}{(2,4,3)}
\psfrag{(3,2,4)}{(3,2,4)}
\psfrag{(4,2,3)}{(4,2,3)}
\psfrag{(3,4,2)}{(3,4,2)}
\psfrag{(4,3,2)}{(4,3,2)}

\psfrag{(12,3,4)}{(12,3,4)}
\psfrag{(12,4,3)}{(12,4,3)}
\psfrag{(13,2,4)}{(13,2,4)}
\psfrag{(13,4,2)}{(13,4,2)}
\psfrag{(14,2,3)}{(14,2,3)}
\psfrag{(14,3,2)}{(14,3,2)}

\psfrag{(23,1,4)}{(23,1,4)}
\psfrag{(23,4,1)}{(23,4,1)}
\psfrag{(24,1,3)}{(24,1,3)}
\psfrag{(24,3,1)}{(24,3,1)}
\psfrag{(34,1,2)}{(34,1,2)}
\psfrag{(34,2,1)}{(34,2,1)}

\psfrag{(4,12,3)}{(4,12,3)}
\psfrag{(3,12,4)}{(3,12,4)}
\psfrag{(4,13,2)}{(4,13,2)}
\psfrag{(2,13,4)}{(2,13,4)}
\psfrag{(3,14,2)}{(3,14,2)}
\psfrag{(2,14,3)}{(2,14,3)}

\psfrag{(4,23,1)}{(4,23,1)}
\psfrag{(1,23,4)}{(1,23,4)}
\psfrag{(3,24,1)}{(3,24,1)}
\psfrag{(1,24,3)}{(1,24,3)}
\psfrag{(2,34,1)}{(2,34,1)}
\psfrag{(1,34,2)}{(1,34,2)}

\psfrag{(3,4,12)}{(3,4,12)}
\psfrag{(4,3,12)}{(4,3,12)}
\psfrag{(2,4,13)}{(2,4,13)}
\psfrag{(4,2,13)}{(4,2,13)}
\psfrag{(2,3,14)}{(2,3,14)}
\psfrag{(3,2,14)}{(3,2,14)}

\psfrag{(1,4,23)}{(1,4,23)}
\psfrag{(4,1,23)}{(4,1,23)}
\psfrag{(1,3,24)}{(1,3,24)}
\psfrag{(3,1,24)}{(3,1,24)}
\psfrag{(1,2,34)}{(1,2,34)}
\psfrag{(2,1,34)}{(2,1,34)}

\includegraphics[width=12cm]{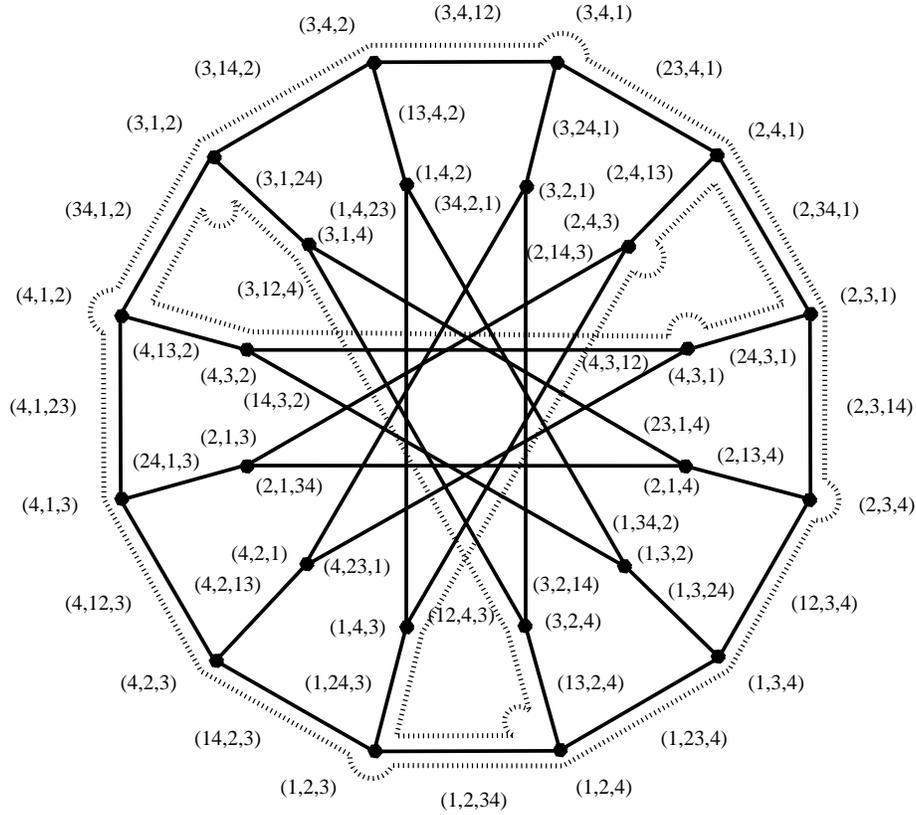}
\end{center}
\caption{The Hom complex ${\rm Hom}(K_3,K_4)$.}
\label{fig:Hom_K3_K4}
\end{figure}

\bigskip

As before in the case of ${\rm Hom}(\overline{C}_5,K_3)$,
every edge of ${\rm Hom}(K_r,K_{r+1})$ can be extended
in exactly two ways to an edge of ${\rm Hom}(\overline{C}_{2r+1},K_{r+1})$.
This can be interpreted geometrically by thickening every edge 
of the $1$-dimensional manifold ${\rm Hom}(K_r,K_{r+1})$
to a $2$-dimensional strip and then gluing these strips together at
the vertices of ${\rm Hom}(K_r,K_{r+1})$. In this way, we get a two-dimensional
manifold with boundary, with the boundary being homeomorphic
to ${\rm Hom}(\overline{C}_{2r+1},K_{r+1})$. 
In Figure~\ref{fig:Hom_K3_K4} we display the Hom complex ${\rm Hom}(K_3,K_4)$,
consisting of $24$ vertices and $36$ edges, together with two of the
$3!=6$ (dotted) cycles  of ${\rm Hom}(\overline{C}_7,K_4)$.
Every vertex of ${\rm Hom}(K_r,K_{r+1})$ can be extended in
$r+1$ ways to an edge of ${\rm Hom}(\overline{C}_{2r+1},K_{r+1})$.
These $r+1$ edges form a path that we display as dotted half-cycles in
the Figures~\ref{fig:C_5_crown} and \ref{fig:Hom_K3_K4}.

\begin{conj} The $3$-dimensional graph coloring manifold\,
${\rm Hom}(\overline{C}_{2r+1},K_{r+2})$, \linebreak
\mbox{$r\geq 2$}, is homeomorphic to the orientable Seifert manifold\, $\{\,Oo,r!-2\,|\,r!\,\}$
with homo\-logy\, $({\mathbb Z},{\mathbb Z}^{2(r!-2)}\oplus{\mathbb Z}_{r!},{\mathbb Z}^{2(r!-2)},{\mathbb Z})$.
\end{conj}

\noindent
The conjecture holds for $r=2$ and $r=3$.
(For an introduction to Seifert manifolds see Seifert \cite{Seifert1933}
as well as \cite{Lutz2003bpre} and \cite{Orlik1972}.)

For $r=2$, Theorem~\ref{thm:rp3} yields $\{\,Oo,0,|\,2\,\}\cong {\mathbb R}{\bf P}^3\cong{{\rm Hom}}(\overline{C}_5,K_4)$.
For $r=3$, the product triangulation of ${\rm Hom}(\overline{C}_7,K_5)$ has $f$-vector $f=(2520,20160,35280,17640)$
and homology $({\mathbb Z},{\mathbb Z}^8\oplus{\mathbb Z}_6,{\mathbb Z}^8,{\mathbb Z})$.
It took ten minutes on a Pentium~R $2.8$ GHz processor to
reduce the triangulation with the \texttt{bistellar} client system of~\cite{polymake} 
to a triangulation
with $f=(27,289,524,262)$. In a second step,
the topological type of the resulting triangulation
was recognized within seconds with the program \texttt{Three-manifold Recognizer}~\cite{Matveev2005} 
(see also \cite{Matveev2003}).
Many thanks to S.~V.~Matveev, E.~Pervova, and V.~Tarkaev
for their help with the recognition!
\begin{thm}
\mbox{}\,${\rm Hom}(\overline{C}_7,K_5)\cong \{\,Oo,4,|\,6\,\}.$
\end{thm}
We will describe further graph coloring manifolds of similar size
in Section~\ref{sec:two}, for which their topological type
was recognized in the same manner.

\medskip

\noindent
\emph{Recognition heuristic for Seifert and graph manifolds:}
\begin{itemize}
\item[1.] Reduce the size of a given triangulation with the \texttt{bistellar} client 
          of the TOPAZ module of the \texttt{polymake} system~\cite{polymake}.
\item[2.] Use the program \texttt{Three-manifold Recognizer}~\cite{Matveev2005} for the recognition.
\end{itemize}

\medskip

\noindent
If the (Matveev) complexity of a given triangulation is not too large,
there is a good chance to recognize the topological type, even
when the triangulation is huge.

\subsection{Hom Complexes of Complements of Even Cycles}

\begin{figure}
\begin{center}
\footnotesize
\psfrag{a}{$a_1$}
\psfrag{b}{$a_2$}
\psfrag{c}{$a_3$}
\psfrag{A}{$A_1$}
\psfrag{B}{$A_2$}
\psfrag{C}{$A_3$}
\includegraphics[height=4.75cm]{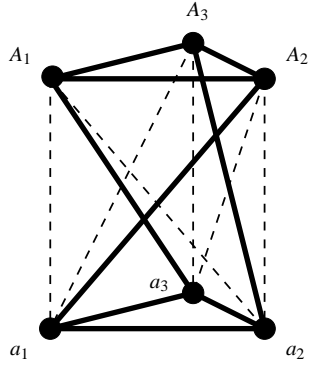}
\end{center}
\caption{The cycle $C_6$ (dashed) and its complement $\overline{C}_6$.}
\label{fig:C_6_C_6_bar}
\end{figure}

Similar to the crown representation of (complements of) odd cycles,
we split the vertices of even cycles $C_{2r}$ into a lower and an upper part,
corresponding to the bipartition of $C_{2r}$. The lower and also the
upper part form a complete graph $K_r$ in $\overline{C}_{2r}$, i.e.,
every maximal cell of ${\rm Hom}(\overline{C}_{2r},K_r)$
contains each color $1,\dots,r$ exactly twice, once in the lower
part and once in the upper part. (Figure~\ref{fig:C_6_C_6_bar} 
displays $C_6$ and its complement $\overline{C}_6$
together with a cell $(a_1,a_2,a_3;A_1,A_2,A_3)$ of ${\rm Hom}(\overline{C}_6,K_r)$.)

We will employ the following two propositions to describe
the $2$-dimensional Hom complexes ${\rm Hom}(\overline{C}_{2r},K_{r+1})$.

\begin{prop} {\rm (Babson and Kozlov~\cite{BabsonKozlov2006})}\label{prp:K_r_K_s}
The Hom complex ${\rm Hom}(K_r,K_s)$
is homotopy equivalent to a wedge 
of $f(r,s)$ spheres of dimension $s-r$,
where the numbers $f(r,s)$ satisfy the recurrence relation
\begin{equation}
f(r,s)=rf(r-1,s-1)+(r-1)f(r,s-1),
\end{equation}
for\, $s>r\geq 2$; with the boundary values
$f(r,r)=r!-1$, $f(1,s)=0$\, for\, $s\geq 1$,
and $f(r,s)=0$\, for\, $r>s$.
\end{prop}

\begin{prop} {\rm (\v{C}uki\'c and Kozlov~\cite{CukicKozlov2005})}
\, $f(r,r+1)=r!\,\frac{r^2-r-2}{2}+1$. 
\end{prop}

\begin{thm}\label{prop:cubical_surface}
The Hom complex ${\rm Hom}(\overline{C}_{2r},K_{r+1})$, $r\geq 2$,
is an orientable cubical surface of genus 
\begin{equation}
\textstyle g(r)=f(r,r+1)=r!\,\frac{r^2-r-2}{2}+1
\end{equation}
with $n(r)=(2+r^2)\cdot(r+1)!$ vertices, $2(n(r)+2g(r)-2)$ edges,
and $n(r)+2g(r)-2$ squares.
\end{thm}

\noindent
\textbf{Proof.}
Let $(a_1,\dots,a_r;A_1,\dots,A_r)$ be a maximal cell of
${\rm Hom}(\overline{C}_{2r},K_{r+1})$.
Since every color $1,\dots,r+1$ appears exactly
once in $(a_1,\dots,a_r)$ and once in $(A_1,\dots,A_r)$
the cell $(a_1,\dots,a_r;A_1,\dots,A_r)$ 
is the product of the edge $(a_1,\dots,a_r)$ 
with the edge $(A_1,\dots,A_r)$. Hence, ${\rm Hom}(\overline{C}_{2r},K_{r+1})$
is a cubical surface. 

We count the vertices of ${\rm Hom}(\overline{C}_{2r},K_{r+1})$.
For every vertex $(v_1,\dots,v_r,w_1,\dots,w_r)$
we have to choose $r$ of the $r+1$ colors for the lower part 
and then have $r!$ choices to place these $r$ colors.
Let $(v_1,\dots,v_r)=(1,\dots,r)$ be such a placement. 
If the left out color $r+1$ does not appear in the upper part, 
then $(1,\dots,r)$ can be extended in exactly two ways to a coloring of $\overline{C}_{2r}$,
yielding the vertices $(1,\dots,r;1,\dots,r)$ and $(1,\dots,r;2,\dots,r,1)$
of ${\rm Hom}(\overline{C}_{2r},K_{r+1})$. If the left out color $r+1$
\emph{is} used in the top part, then there are $r$ choices to place it,
and for each such placement every choice to not use one of
the colors $1,\dots,r$ determines a vertex. Therefore, we have 
altogether $2+r^2$ choices 
to extend $(1,\dots,r)$ to a vertex of ${\rm Hom}(\overline{C}_{2r},K_{r+1})$;
i.e., ${\rm Hom}(\overline{C}_{2r},K_{r+1})$
has\, $n(r):=(2+r^2)\cdot(r+1)!$ vertices.

Let $M$ be an orientable cubical surface of genus $g$
with $n$ vertices, $e$ edges, and $s$ squares.
Since every square is bounded by four edges
and every edge appears in two squares, 
double counting yields $2e=4s$.
By this equation and by Euler's relation, $s-e+n=\chi(M)=2-2g$,
we get that $s=n+2g-2$ and $e=2(n+2g-2)$.

It remains to show that ${\rm Hom}(\overline{C}_{2r},K_{r+1})$
is orientable and has genus $g(r)=f(r,r+1)=r!\,\frac{r^2-r-2}{2}+1$.
For this, let us fix an edge, say $(a_1,\dots,a_r)=(1,2,\dots,r-1,r(r+1))$,\linebreak
of ${\rm Hom}(K_r,K_{r+1})$.
Then the sequence of $2r$ squares

\bigskip

\begin{tabular}{l}
$(1,2,\dots,r-1,r(r+1);1,2,\dots,r-2,r-1,r(r+1))$,\\
$(1,2,\dots,r-1,r(r+1);1,2,\dots,r-2,(r-1)r,r+1)$,\\
$(1,2,\dots,r-1,r(r+1);1,2,\dots,(r-2)(r-1),r,r+1)$,\\
\dots\\
$(1,2,\dots,r-1,r(r+1);1,23,\dots,r-1,r,r+1)$,\\
$(1,2,\dots,r-1,r(r+1);12,3,\dots,r-1,r,r+1)$,\\
$(1,2,\dots,r-1,r(r+1);2,3,\dots,r-1,r,1(r+1))$,\\
$(1,2,\dots,r-1,r(r+1);2,3,\dots,r-1,r(r+1),1)$,\\
$(1,2,\dots,r-1,r(r+1);2,3,\dots,r-1,r+1,1r)$,\\
$(1,2,\dots,r-1,r(r+1);12,3,\dots,r-1,r+1,r)$,\\
$(1,2,\dots,r-1,r(r+1);1,23,\dots,r-1,r+1,r)$,\\
\dots\\
$(1,2,\dots,r-1,r(r+1);1,2,\dots,(r-2)(r-1),r+1,r)$,\\
$(1,2,\dots,r-1,r(r+1);1,2,\dots,r-2,(r-1)(r+1),r)$\\
\end{tabular}

\begin{figure}
\begin{center}
\footnotesize
\includegraphics[height=6.5cm]{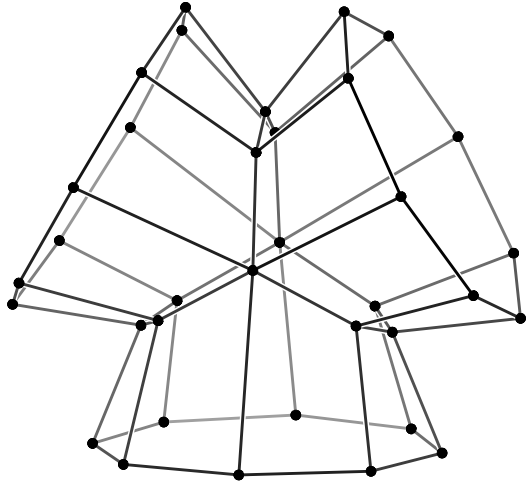}
\end{center}
\caption{Three cylinders forming a trinoid in ${\rm Hom}(\overline{C}_6,K_4)$.}
\label{fig:C_5_C_5_bar}
\end{figure}

\bigskip

\noindent
forms a cylinder $C_{2r}\times I$. By symmetry, we get such a
cylinder for every edge $(a_1,\dots,a_r)$ of ${\rm Hom}(K_r,K_{r+1})$.
Since every vertex of the graph ${\rm Hom}(K_r,K_{r+1})$
has degree $r$, we have $r$ cylinders in ${\rm Hom}(\overline{C}_{2r},K_{r+1})$
meeting ``at a vertex'' of ${\rm Hom}(K_r,K_{r+1})$. (In the case of ${\rm Hom}(K_3,K_4)$ three cylinders meet
at a vertex, which yields a trinoid as depicted in Figure~\ref{fig:C_5_C_5_bar}.)
By inspecting the gluing at the vertices, it is easy to deduce that
${\rm Hom}(\overline{C}_{2r},K_{r+1})$ is orientable.
It moreover follows that ${\rm Hom}(\overline{C}_{2r},K_{r+1})$ has genus $f(r,r+1)$,
which is the number of wedged $1$-spheres in the graph ${\rm Hom}(K_r,K_{r+1})$.
\mbox{}\hfill$\Box$

\bigskip

As in the case of ${\rm Hom}(\overline{C}_{2r+1},K_{r+1})$,
we can interpret ${\rm Hom}(\overline{C}_{2r},K_{r+1})$ geometrically
in the following way. If we thicken the edges of the $1$-dimensional
manifold ${\rm Hom}(K_r,K_{r+1})$ to solid tubes, 
then for the resulting $3$-manifold with boundary
the boundary is homeomorphic to ${\rm Hom}(\overline{C}_{2r},K_{r+1})$.

\begin{conj}\label{conj:connected_sum}
The Hom complex\, ${\rm Hom}(\overline{C}_{2r},K_s)$\,
is, for $s>r\geq 2$, homeomorphic to the connected sum of\,\, $f(r,s)$\, copies
of\,\, $S^{s-r}\!\times S^{s-r}$.
\end{conj}

The $4$-dimensional Hom complex ${{\rm Hom}}(\overline{C}_6,K_5)$ consists of $3180$ cells
and has $1920$ vertices. The corresponding product triangulation
has $f=(1920,30780,104520,126000,$ $50400)$.
With the \texttt{bistellar} client it took half a day to
reduce this triangulation to a triangulation with
$f$-vector $(33,379,1786,2300,920)$. 
The latter triangulation is simply connected, as we computed with the 
group algebra package \texttt{GAP} \cite{GAP4}.
The homology of the triangulation is
$({\mathbb Z},0,{\mathbb Z}^{58},0,{\mathbb Z})$.
Moreover, we used \texttt{polymake} to compute the intersection form
of the example, which turned out to be indefinite, even, and of rank $29$.
By the classification of Freedman~\cite{Freedman1982},
this shows:

\begin{thm} ${\rm Hom}(\overline{C}_6,K_5)\cong (S^2\!\times S^2)^{\# 29}$.
\end{thm}

Conjecture~\ref{conj:connected_sum} thus
holds in the case $(r,s)=(3,5)$, and,
by Theorem~\ref{prop:cubical_surface},
also for the series $s=r+1\geq 3$.

\section{Graph Coloring Manifolds of Sphere Dimension Two}
\label{sec:two}

Flag simplicial $2$-spheres with small numbers $n$ of vertices can be obtained
by first enumerating \emph{all} triangulated $2$-spheres with $n$
vertices and then testing which of these are flag.
Triangulations of two-dimensional spheres with up to
$23$ vertices have been enumerated with the program
\texttt{plantri} by Brinkmann and McKay~\cite{plantri} 
(see the manual of \texttt{plantri} or the web-page of Royle~\cite{Royle_url}
for the numbers of triangulations on $n\leq 23$ vertices).
With another approach, triangulations of all two-dimensional manifolds with up to $10$ vertices 
have been enumerated by Lutz (cf.~\cite{Lutz2005apre});
the respective numbers of triangulations are given in Table~\ref{tbl:number_of_triangulations}.

\begin{figure}[!ht]
\begin{center}
\footnotesize
\includegraphics[height=2.75cm]{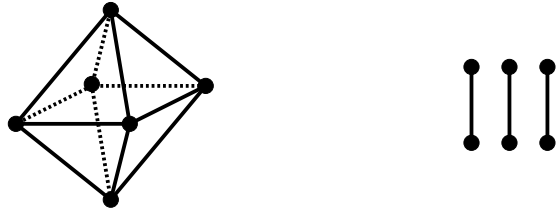}
\end{center}
\caption{The flag sphere $\mbox{}^2\hspace{.3pt}6_{\,3}=\partial C_3^{\Delta}$ and the complement of its $1$-skeleton.}
\label{fig:n6_flag_black}
\end{figure}

\begin{figure}[!ht]
\begin{center}
\footnotesize
\includegraphics[height=2.75cm]{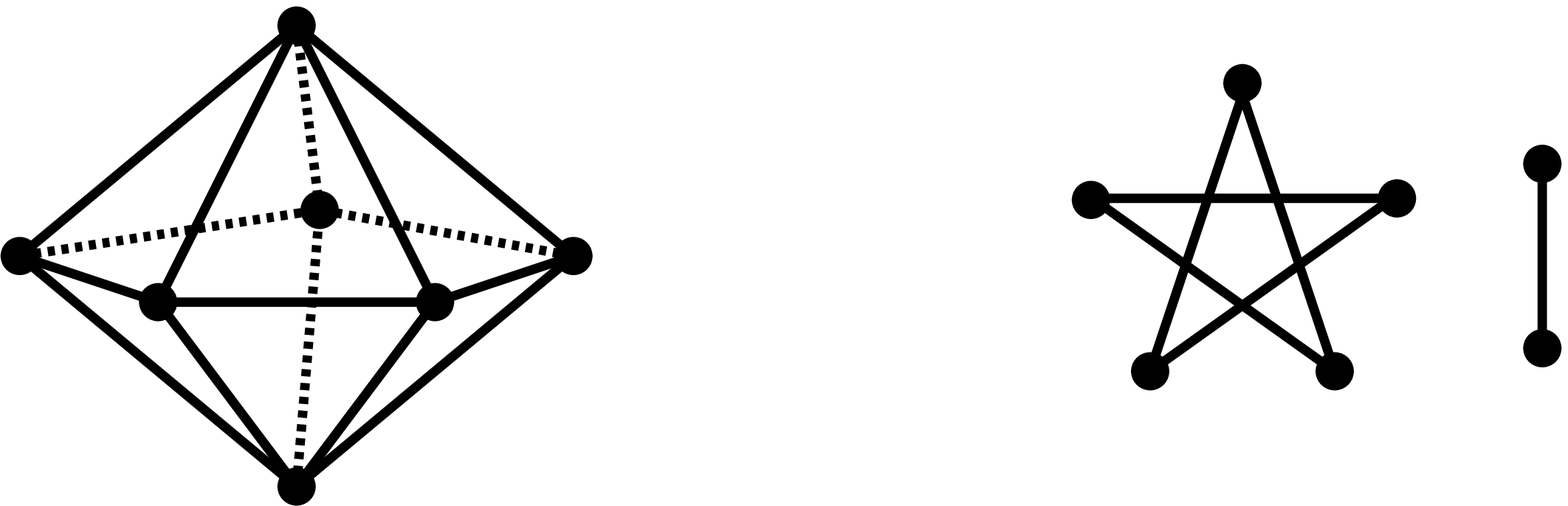}
\end{center}
\caption{The flag sphere $\mbox{}^2\hspace{.3pt}7_{\,9}=C_5*S^0$ and the complement of its $1$-skeleton.}
\label{fig:n7_flag_black}
\end{figure}

\begin{figure}
\begin{center}
\footnotesize
\includegraphics[height=2.75cm]{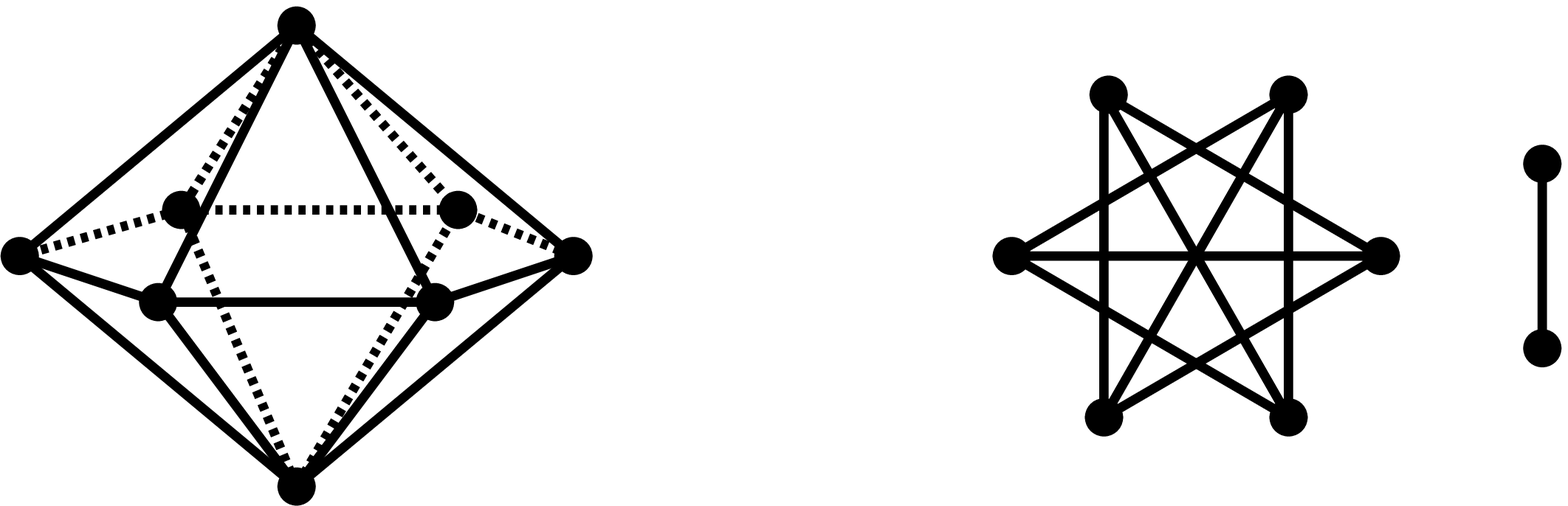}
\end{center}
\caption{The flag sphere $\mbox{}^2\hspace{.3pt}8_{\,41}=C_6*S^0$ and the complement of its $1$-skeleton.}
\label{fig:n8_41_flag_black}
\end{figure}

\begin{figure}
\begin{center}
\footnotesize
\includegraphics[height=2.75cm]{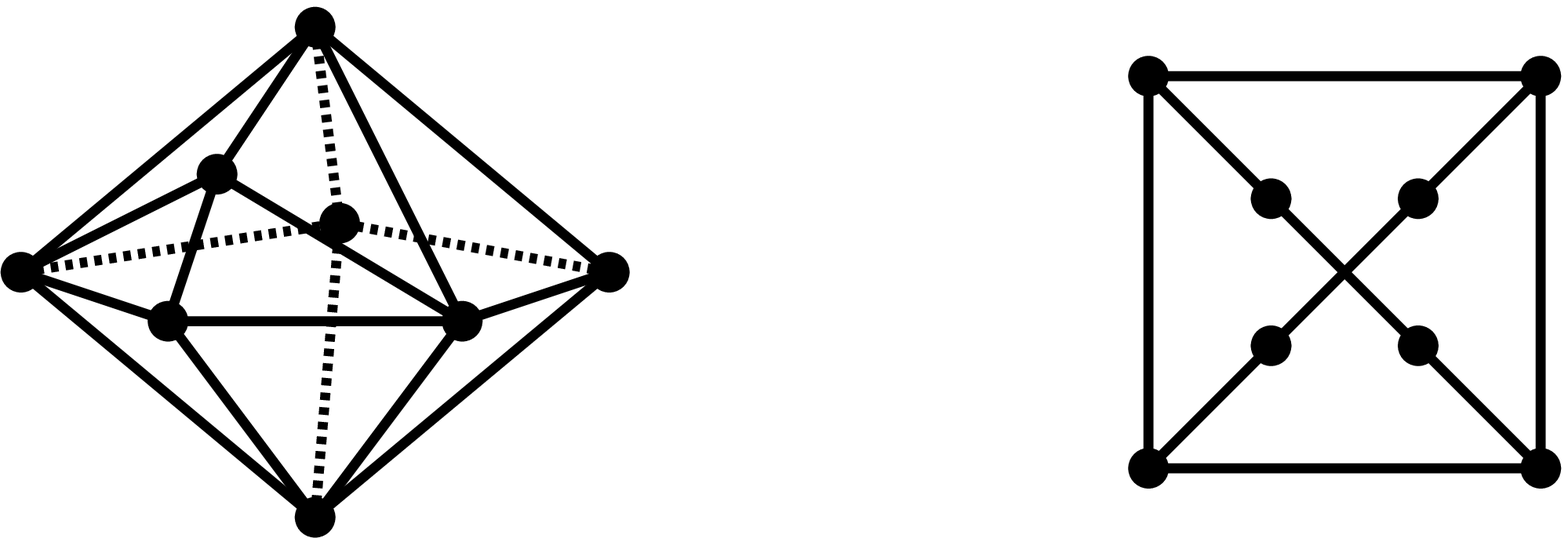}
\end{center}
\caption{The flag sphere $\mbox{}^2\hspace{.3pt}8_{\,43}$ and the complement of its $1$-skeleton.}
\label{fig:n8_43_flag_black}
\end{figure}

\begin{figure}
\begin{center}
\footnotesize
\includegraphics[height=2.75cm]{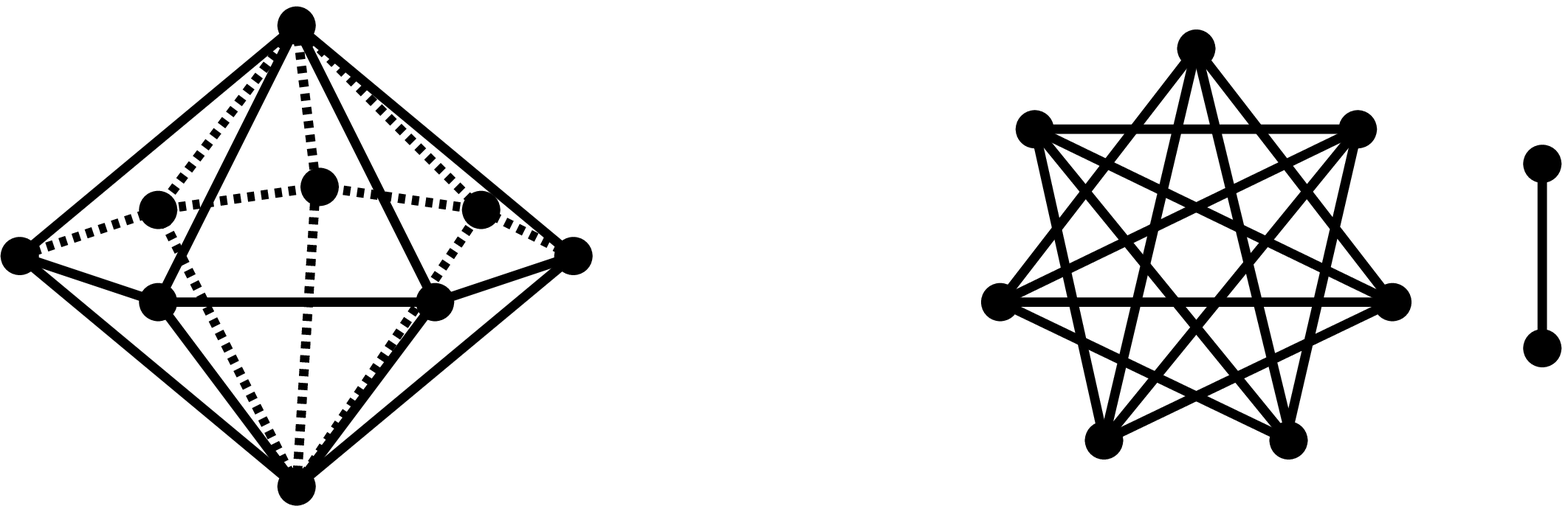}
\end{center}
\caption{The flag sphere $\mbox{}^2\hspace{.3pt}9_{\,630}=C_7*S^0$ and the complement of its $1$-skeleton.}
\label{fig:n9_630_flag_black}
\end{figure}

\begin{figure}
\begin{center}
\footnotesize
\includegraphics[height=2.75cm]{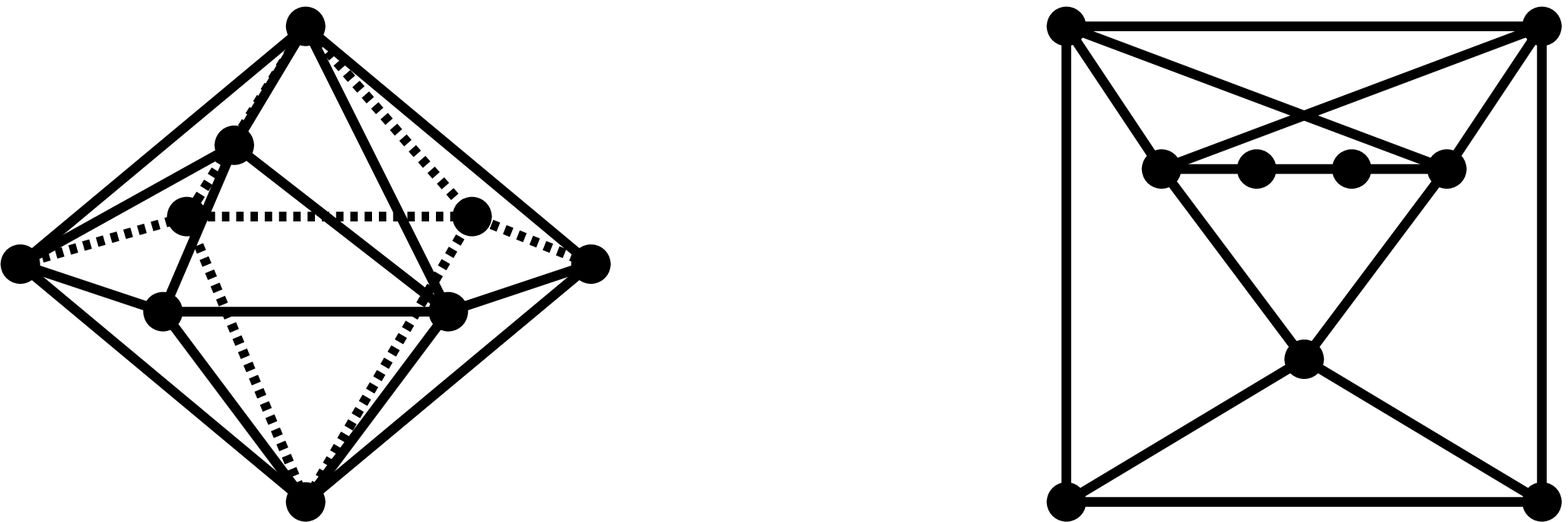}
\end{center}
\caption{The flag sphere $\mbox{}^2\hspace{.3pt}9_{\,651}$ and the complement of its $1$-skeleton.}
\label{fig:n9_651_flag_black}
\end{figure}

\begin{figure}
\begin{center}
\footnotesize
\includegraphics[height=2.75cm]{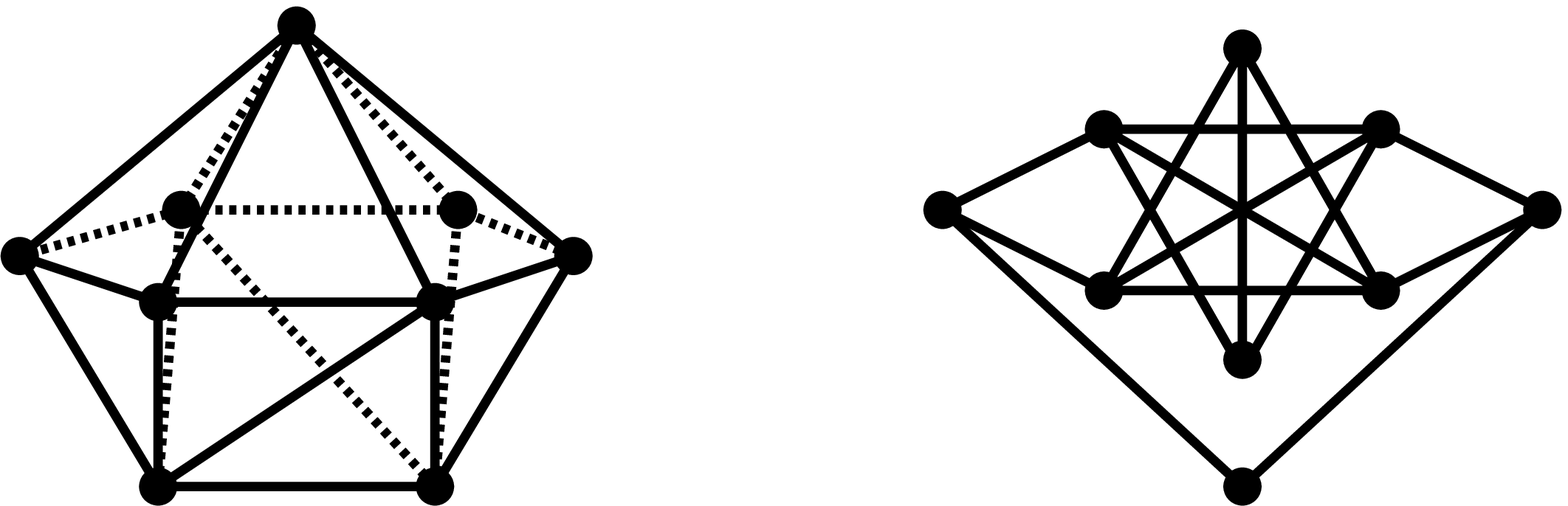}
\end{center}
\caption{The flag sphere $\mbox{}^2\hspace{.3pt}9_{\,652}$ and the complement of its $1$-skeleton.}
\label{fig:n9_652_flag_black}
\end{figure}

\begin{figure}
\begin{center}
\footnotesize
\includegraphics[height=3.7cm]{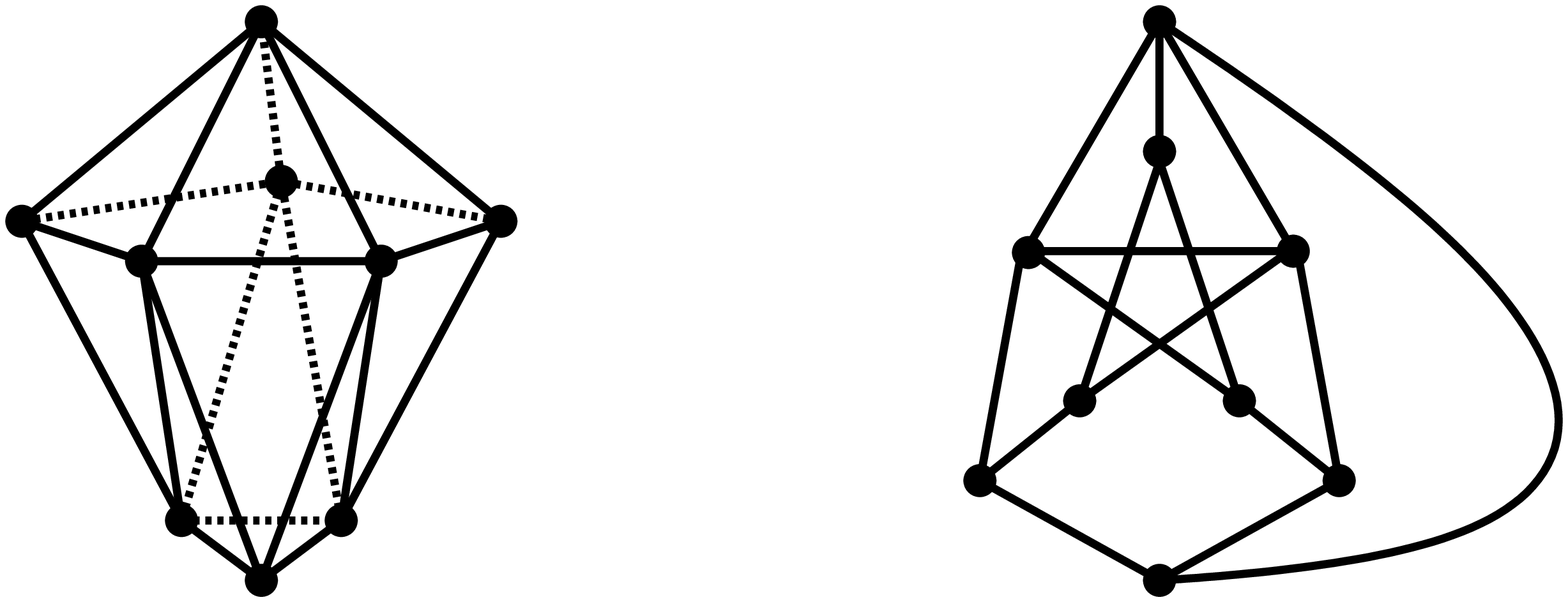}
\end{center}
\caption{The flag sphere $\mbox{}^2\hspace{.3pt}9_{\,652}$ and the complement of its $1$-skeleton.}
\label{fig:n9_655_flag_black}
\end{figure}

\begin{table}
\small\centering
\defaultaddspace=0.3em
\caption{Triangulated surfaces with few vertices.}\label{tbl:number_of_triangulations}
\begin{tabular}{l@{\hspace{.75cm}}r@{\hspace{.75cm}}r@{\hspace{.75cm}}r@{\hspace{.75cm}}r@{\hspace{.75cm}}r}
 \addlinespace
 \addlinespace
 \addlinespace
 \addlinespace
\toprule
 \addlinespace
\# Vertices      & 6 & 7 &  8 &   9 &    10 \\
\midrule
 \addlinespace
 \addlinespace
 \addlinespace
\# Manifolds     & 3 & 9 & 43 & 655 & 42426 \\
 \addlinespace
\# Spheres       & 2 & 5 & 14 &  50 &   233 \\
 \addlinespace
\# Flag Spheres  & 1 & 1 &  2 &   4 &    10 \\
 \addlinespace
 \addlinespace
\bottomrule
\end{tabular}
\end{table}
The flag simplicial spheres with up to $9$ vertices
together with the complements of their $1$-skeleta
are displayed in Figures~\ref{fig:n6_flag_black}--\ref{fig:n9_655_flag_black}.
(The symbol $\mbox{}^2\hspace{.3pt}n_{\,k}$ stands for the
$k$th $2$-manifold with $n$ vertices in the catalog~\cite{Lutz_PAGE}.)

For the flag $2$-spheres $\mbox{}^2\hspace{.3pt}6_{\,3}=\partial C_3^{\Delta}$
(the boundary of the $3$-dimensional cross-poly\-tope~$C_3^{\Delta}$),
$\mbox{}^2\hspace{.3pt}7_{\,9}=C_5*S^0$,
$\mbox{}^2\hspace{.3pt}8_{\,41}=C_6*S^0$, and
$\mbox{}^2\hspace{.3pt}9_{\,630}=C_7*S^0$,
the complements of the respective $1$-skeleta are not connected,
and therefore, by Equation~\ref{eq:direct_product}, are direct products.

For those flag $2$-spheres with $n\leq 9$ vertices,
for which the complements of their $1$-skeleta
are connected, we analyzed the product triangulations 
of their Hom complexes with few colors. 
Table~2  gives the results.

\begin{table}\label{tbl:flagtwo}.
\small\centering
\defaultaddspace=0.3em
\caption{Hom complexes associated with the flag spheres
  $\mbox{}^2\hspace{.3pt}8_{\,43}$,
  $\mbox{}^2\hspace{.3pt}9_{\,651}$,
  $\mbox{}^2\hspace{.3pt}9_{\,652}$, and $\mbox{}^2\hspace{.3pt}9_{\,655}$.}
\begin{tabular*}{\linewidth}{@{\extracolsep{\fill}}llll@{}}
 \addlinespace
 \addlinespace
 \addlinespace
 \addlinespace
\toprule
 \addlinespace
    Hom complex         & Type & Homology & $f$-Vector of Product Triangulation \\
\midrule
 \addlinespace
 \addlinespace
 \addlinespace
${\rm Hom}(\overline{{\rm SK}_1(\mbox{}^2\hspace{.3pt}8_{\,43})},K_3)$ & $4$ cycles   &  & each cycle has $24$ vertices \\
 \addlinespace
${\rm Hom}(\overline{{\rm SK}_1(\mbox{}^2\hspace{.3pt}9_{\,651})},K_3)$ & $24$ vertices &  & \\
 \addlinespace
${\rm Hom}(\overline{{\rm SK}_1(\mbox{}^2\hspace{.3pt}9_{\,652})},K_3)$ & $24$ vertices &  & \\
 \addlinespace
${\rm Hom}(\overline{{\rm SK}_1(\mbox{}^2\hspace{.3pt}9_{\,655})},K_3)$ & $12$ vertices &  & \\
 \addlinespace
 \addlinespace
\midrule
 \addlinespace
 \addlinespace
 \addlinespace
${\rm Hom}(\overline{{\rm SK}_1(\mbox{}^2\hspace{.3pt}8_{\,43})},K_4)$  & ?                                   & $({\mathbb Z},{\mathbb Z}\oplus{\mathbb Z}_2,{\mathbb Z}_2,{\mathbb Z},{\mathbb Z})$ & $(3624,55224,184656,221760,88704)$ \\
 \addlinespace
${\rm Hom}(\overline{{\rm SK}_1(\mbox{}^2\hspace{.3pt}9_{\,651})},K_4)$ & $({\mathbb T}^2)^{\# 13}\!\times\!S^1$ & $({\mathbb Z},{\mathbb Z}^{27},{\mathbb Z}^{27},{\mathbb Z})$                        & $(2928,21360,36864,18432)$ \\
 \addlinespace
${\rm Hom}(\overline{{\rm SK}_1(\mbox{}^2\hspace{.3pt}9_{\,652})},K_4)$ & $({\mathbb T}^2)^{\# 13}\!\times\!S^1$ & $({\mathbb Z},{\mathbb Z}^{27},{\mathbb Z}^{27},{\mathbb Z})$                        & $(3120,22992,39744,19872)$ \\
 \addlinespace
${\rm Hom}(\overline{{\rm SK}_1(\mbox{}^2\hspace{.3pt}9_{\,655})},K_4)$ & $(S^2\!\times\!S^1)^{\# 13}$           & $({\mathbb Z},{\mathbb Z}^{13},{\mathbb Z}^{13},{\mathbb Z})$                        & $(3096,22104,38016,19008)$ \\
 \addlinespace
 \addlinespace
\bottomrule
\end{tabular*}
\end{table}

\smallskip

\subsection*{Acknowledgements} 

The authors are grateful to S.~Felsner for helpful discussions. 
Many thanks also to S.~V.~Mat\-veev, E.~Pervova, and V.~Tarkaev
for their help with the recognition of\, $3$-dimensional graph coloring manifolds.
Moreover, we thank the anonymous referee for helpful remarks
that led to a substantial improvement of the display of Section~\ref{sec:main}.

\smallskip

\subsection*{\normalsize{Note added in proof}} 
Conjecture~\ref{conj:csorba} has recently been proved by C.~Schultz~\cite{Schultz2005bpre}.

\medskip


\end{document}